\documentclass{article}                     
\pdfoutput=1
\usepackage{arxiv}
\usepackage{graphicx}
\usepackage{hyperref}
\usepackage{natbib}
\usepackage[utf8]{inputenc}
\usepackage{amsmath}  
\usepackage{xcolor} 
\usepackage{amssymb} 
\usepackage{bm} 
\usepackage{tikz} 
\usetikzlibrary{positioning}
\usepackage{subcaption} 
\usepackage{algorithmicx} 
\usepackage{algpseudocode}
\usepackage[parfill]{parskip}

\newcommand{\norm}[1]{\left\lVert#1\right\rVert}

\newcommand{\ybm}{\bm{y}}

\newcommand{\xbm}{\bm{x}}

\newcommand{\wbm}{\bm{w}}
\newcommand{\Wbm}{\bm{W}}
\newcommand{\ubm}{\bm{u}}

\newcommand{\ybmh}{\hat{\bm{y}}}

\newcommand{\ubmh}{\hat{\bm{u}}}

\newcommand{\Qbm}{\bm{Q}}

\newcommand{\Ibm}{\bm{I}}
\newcommand{\Bbm}{\bm{B}}
\newcommand{\bbm}{\bm{b}}

\newcommand{\fbm}{\bm{f}}

\newcommand{\Sigmabm}{\bm{\Sigma}}

\newcommand{\xibm}{\bm{\xi}}

\newcommand{\IR}{\mathbb{R}}

\title{Weighted high dimensional data reduction of
  finite Element Features - An Application on High Pressure of an Abdominal Aortic Aneurysm}

\date{May 1, 2023}

\author{ \href{https://orcid.org/0000-0003-4313-3623}{\includegraphics[scale=0.06]{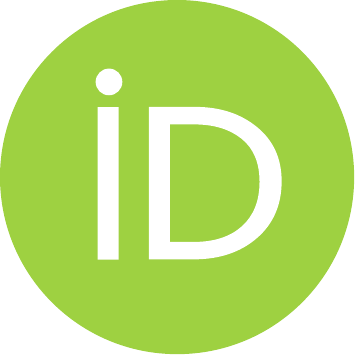}\hspace{1mm}Christoph Striegel}\\
  Department of Statistics\\
  Faculty of Mathematics, Informatics and Statistics\\
  Ludwig-Maximilians-Universität München\\
  Ludwigstr. 33, 80539 Munich\\
	\texttt{christoph.striegel@stat.uni-muenchen.de} \\
	\And
	Göran Kauermann \\
  Department of Statistics\\
  Faculty of Mathematics, Informatics and Statistics\\
  Ludwig-Maximilians-Universität München\\
  Ludwigstr. 33, 80539 Munich\\
	\AND
  Jonas Biehler\\
  Institute for Computational Mechanics\\
  Department of Mechanical Engineering\\
  Technical University of Munich\\
  Munich, Germany
}

\renewcommand{\shorttitle}{Finite element features using additional weights}

\begin{document}
\maketitle

\begin{abstract}
  In this work we propose a low rank approximation of high fidelity finite element simulations
  by utilizing weights corresponding to areas of high stress levels for an abdominal aortic aneurysm, i.e. a
  deformed blood vessel.
  We focus on the van Mises stress, which corresponds to the rupture risk of the aorta.
  This is modeled as a Gaussian Markov random field and we
  define our approximation as a basis of vectors that solve a series of optimization problems.
  Each of these problems describes the minimization of an expected weighted quadratic loss.
  The weights, which encapsulate the importance of each grid point of the finite elements,
  can be chosen freely - either data driven or by incorporating domain knowledge.
  Along with a more general discussion of mathematical properties we provide an effective numerical
  heuristic to compute the basis under general conditions.
  We explicitly explore two such bases on the surface of a high fidelity finite element grid and show
  their efficiency for compression.
  We further utilize the approach to predict the van Mises stress in areas of interest
  using low and high fidelity simulations.
  Due to the high dimension of the data we have to take extra care to keep the problem numerically feasible.
  This is also a major concern of this work.
\end{abstract}


\section{Introduction}
The starting point for the method which we develop in this work is
a real world application whose nature is tightly connected to
the concrete implementation. To be specific, we deal with data
from a computationally expensive high fidelity computer simulation
that returns the deformations and stresses of an abdominal aortic aneurysm
(AAA) in response to blood pressure, see \cite{Biehler2014}.
This simulation is complex, yet deterministic and returns the vector
$\ybm=\left(\ybm(s_{1}),...,\ybm(s_{m_y}) \right)^T, \quad \ybm \in \IR^{m_y}$,
with $\ybm(s_{j})$ denoting the outcome at location $s_{j} \in \mathcal{S}_{y}$,
where $\mathcal{S}_{y}$ is a three-dimensional finite element grid.
The outcome quantity is thereby the magnitude of the van Mises stress,
which is a mechanical quantity that summarizes the tensor-valued stress at each
grid point into one number.
The simulation depends on a (tuning) parameter vector of wall parameters $\xibm$ and we may therefore
denote the outcome of the simulation as $\ybm_{\xibm}$.
\begin{figure}[h]
  \centering
  \begin{subfigure}[b]{1\textwidth}
    \centering
    \includegraphics[width=1.05\textwidth]{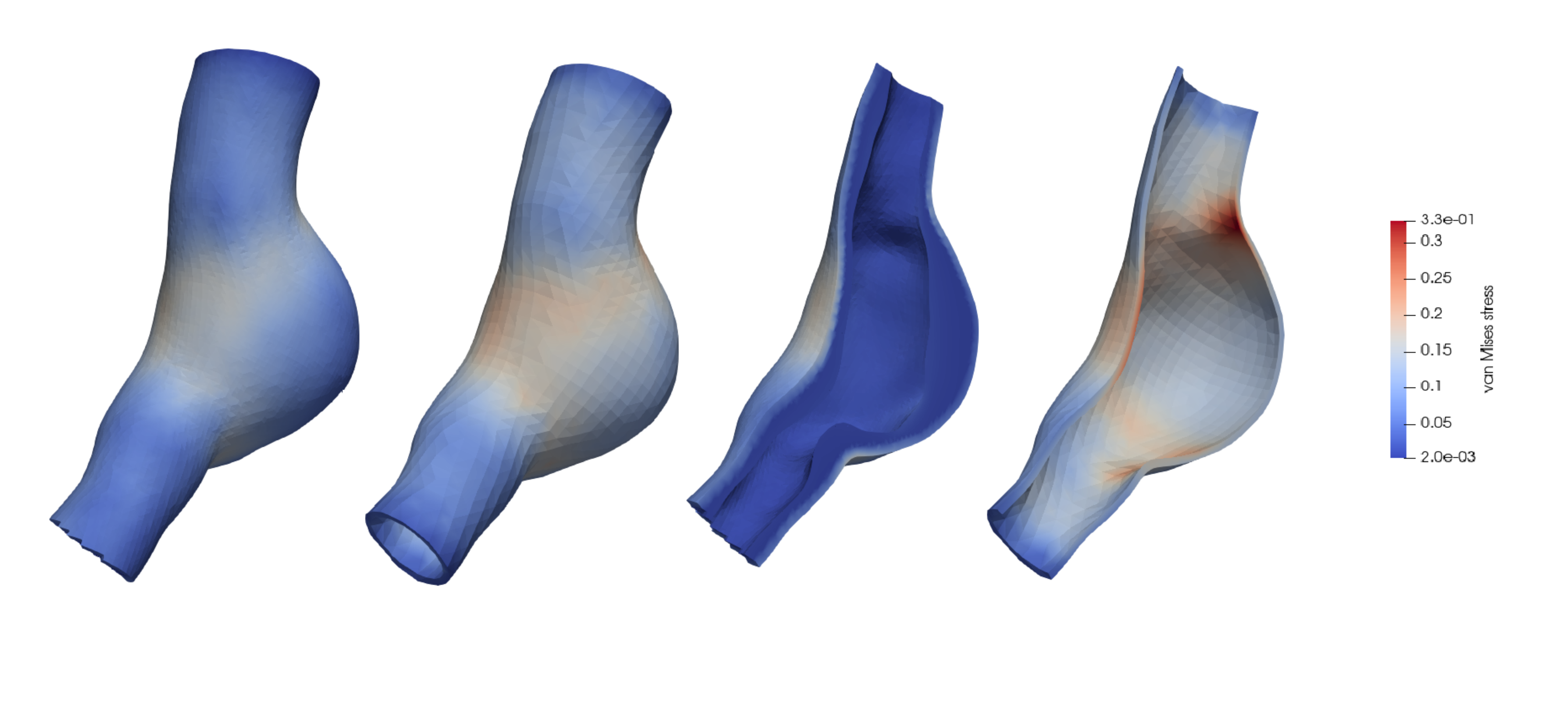}
  \end{subfigure}
  \caption{High (first and third image) and low (second and fourth image) fidelity output for the same set of parameters.
    The images show the artery from the outside (two images to the left) and cutaway
    (two images to the right), demonstrating the three-dimensional structure of the
    data. In addition to other simplifications the low fidelity model is also
    computed using a coarser mesh (not shown above) than the high fidelity model.}
  \label{fig: low_and_high_fid_ex.pdf}
\end{figure}
The outcome of a single simulation is shown as part of Figure \ref{fig: low_and_high_fid_ex.pdf} (first and third image).
These finite element simulations can be used to assess the rapture risk of the blood vessel and to
guide medical treatment in a real word scenario.

Because the high fidelity simulation is computationally expensive
a two step approach has been proposed in \cite{Striegel2022} (see also \cite{Biehler2016})
that uses adjusted low fidelity simulations on a simplified
grid as a surrogate for the high fidelity computations.
The low fidelity simulation can also be seen as part of Figure \ref{fig: low_and_high_fid_ex.pdf} (second and fourth image).
We denote the low fidelity simulation with $\xbm(\xibm)$.
While it can be useful to relate the entire low fidelity simulation
$\xbm$ to its high fidelity counterpart $\ybm$, as done in
\cite{Striegel2022}, it is omitted that one is mostly
interested in high stress levels, as those indicate the risk
of an aortic rupture,
leading to internal bleeding, which, if not treated immediately,
results in the death of the patient.
In other words one aims not to predict the entire high fidelity
outcomes $\ybm_{\xibm}$ but only
high stress areas based on a low fidelity simulation $\xbm_{\xibm}$.
This is the goal of the following paper.
To do so we pursue a weighted dimension reduction of $\ybm_{\xibm}$
through an optimally chosen basis, where the weights account for
relevancy of each grid point
on the high fidelity mesh
with the focus on high stress and hence high risk of rupture.
The basis is computed in an iterative way, where each basis vector
solves an optimization problem based on the weighted expected quadratic compression error.
While the unweighted problem is solved by simple eigenvectors, the addition of
weights makes a new numerical solution necessary.
As each basis vector has the dimension of the high fidelity grid,
i.e. multiple thousand, the
optimization is tedious and the numerical solution, which is highly
optimized to our problem at hand,
is a major contribution of this paper.
However, the theoretical derivations given in the paper are applicable
to any high dimensional Gaussian random field and even more broadly to
any field with a given covariance structure.

Conceptually related work to our basis computation can be found in \cite{Ruben1979}, who present the
approximation of a matrix with additional weights using a weighted least square loss. \cite{Tamuz2005}
uses this approach to develop an effective algorithm to filter out
systematic effects in light curves.
Similarly to both papers, our work can be seen as a generalization of principal component analysis and uses
numerical optimization instead of closed form solutions.
Our methodology does assume a Gaussian random field and not a matrix as starting point
and also our loss function incorporates the elementwise error for every component of the vector - not
the cumulative error for all matrix elements.
In this line, we consider our work as a straight forward extension to (general) eigenvalue decomposition and similar methods
as presented in \cite{ghojogh2019eigenvalue} and similarly it
allows for dual formulations - though not as
a matrix but as a more involved tensor decomposition.
Taking the view of popular fields like computer vision and machine learning, the basis computation
can be seen as the extraction of features from the grid with respect to a certain loss function.
See \cite{zebari2020} for a comprehensive survey or \cite{Soltanpour2017} for the more related usecase of
feature extraction from 3D face images.
In contrast to most methods applied in these fields however our work
is designed to deal with very small sample sizes
in the order of hundreds and even less, where machine
learning typically is not applicable.

The remainder of this work is structured as follows:
First, we will justify our modeling approach by discussing the
nature of the data at hand.
In the third part we will derive our data driven weights model which will
deliver additional information for the weighted basis
approach discussed in the fourth section. Here we define relevant problems
for our basis vectors and provide ideas for an efficient solution along with
exemplary basis vectors for our application at hand and an evaluation of their
quality.
We also touch the numerics behind the basis computation.
Finally, in the fifth section, we evaluate the use of the weighted basis vectors
for an exemplary application, the prediction of high fidelity stress values from low fidelity
counterparts.
The paper finishes with a brief discussion.





\section{Review of the data}
\label{sec:Nature of the maximum}
We define with $\ybm_{<k>}$ the $k$-th high fidelity simulation
run with input parameter $\xibm_{k}$
and accordingly with $\xbm_{<k>}$ the matching low
fidelity counterpart.
The high fidelity outcome $\ybm$ is of dimension $80000$ and
the dimension of the low fidelity counterpart $\xbm$ is $3500$.
Note that the computation of a single high fidelity outcome takes
approximately $10800$ CPU seconds (around $3$ hours),
whereas the computation of a low fidelity simulation takes a mere $215$
CPU seconds (around $3.5$ minutes). This readily shows the
value of a model that enables high fidelity predictions from low fidelity outcomes.
The ultimate goal of the developed methodology is therefore to work based on
small samples with only a limited amount of high and low fidelity
simulations available.
To be specific, we will
only use few cases for the training and define the training set
as $\mathcal{N}_{train}$, with $n_{train} = |\mathcal{N}_{train}| = 50$.
The data at hand trace from a larger experiment and we
utilize the remaining data, which consists of $850$ observations
for evaluating the performance of our approach.
We denote the test set with $\mathcal{N}_{test}$.

As outlined in the introduction we are interested in high
(van Mises) stress values.
In what follows we will exemplary look at the top $0.1\%$ and $0.5\%$ values.
In order to analyze the location where these values occur
we plot the simple high stress frequency, i.e. the number of
samples in the test set for which a high stress value is taken at a certain point, normalized by the number
of simulations:
\begin{equation}
\bm{f}(s_{j})=|k\in \mathcal{N}_{train}:
\ybm_{<k>}(s_{j})>q|/n_{train},
\label{eq: stress frequency}
\end{equation}
for a given quantile $q$ and a point $s_{j}$.
The results for both percentiles of interest are shown
in Figure \ref{fig: location high stress values}:

\begin{figure}[h!]
  \begin{center}
    \includegraphics[width=0.75\linewidth]{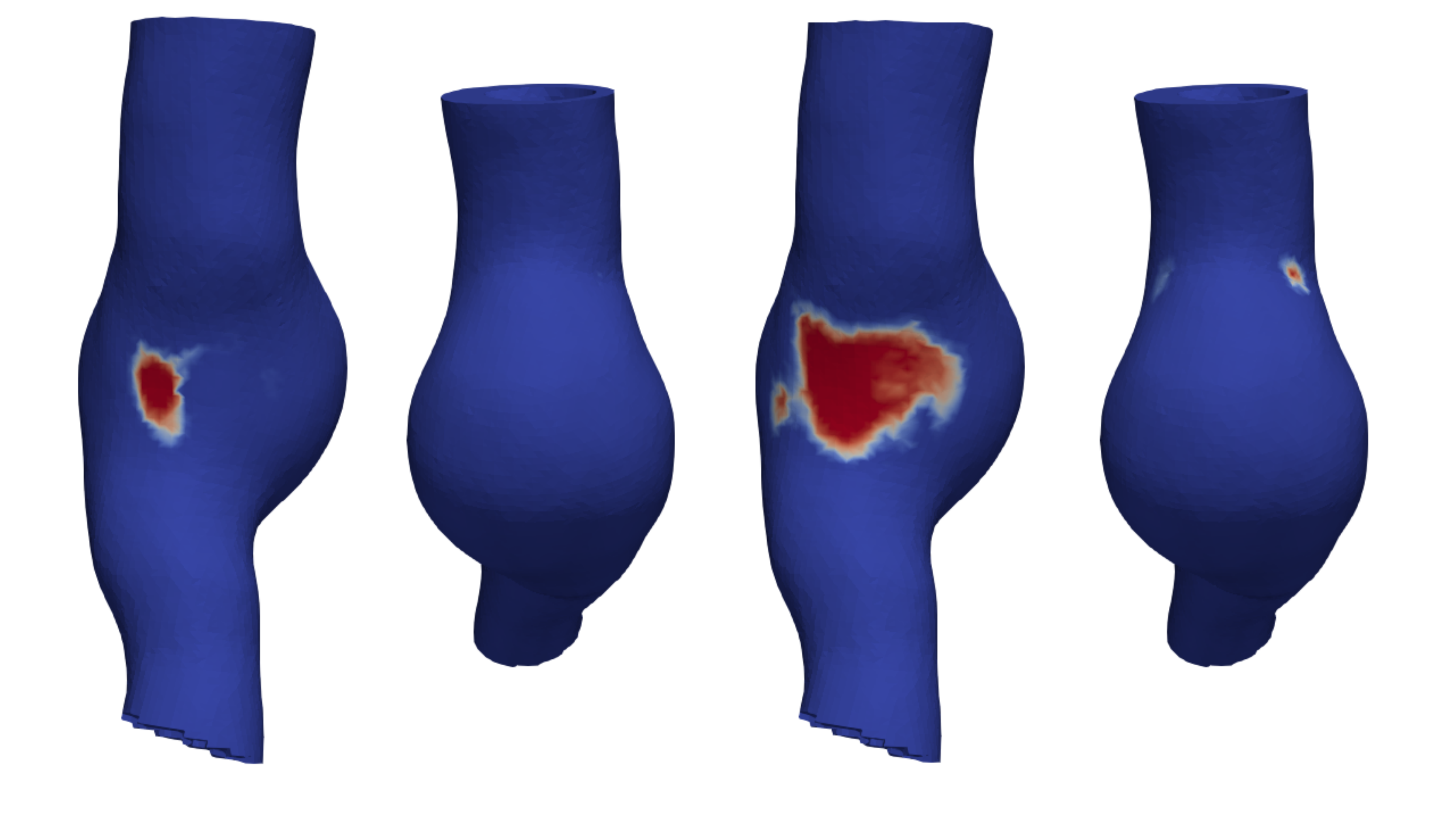}
  \end{center}
  \vspace{-10pt}
  \caption{Average location
    frequency of the top $0.1\%$ (first and second) and top $0.5\%$ (third and fourth) stress values.
  Both front and back view.}
  \label{fig: location high stress values}
\end{figure}

We can see that in both cases the values are spatially clustered on the
grid. In particular there are no isolated grid points.
While in the case of the top $0.1\%$
stress values only the front part of the geometry is of interest, for the top $0.5\%$ there are
additional smaller stress centers at the back.
The fact that the area at the crook has significantly
lower stress may be surprising at first glance, but is in fact a direct consequence of the thrombus,
i.e. a layer of coagulated blood that strengthens this part of the vessel.

\begin{figure}[h!]
  \begin{center}
    \includegraphics[width=0.6\linewidth]{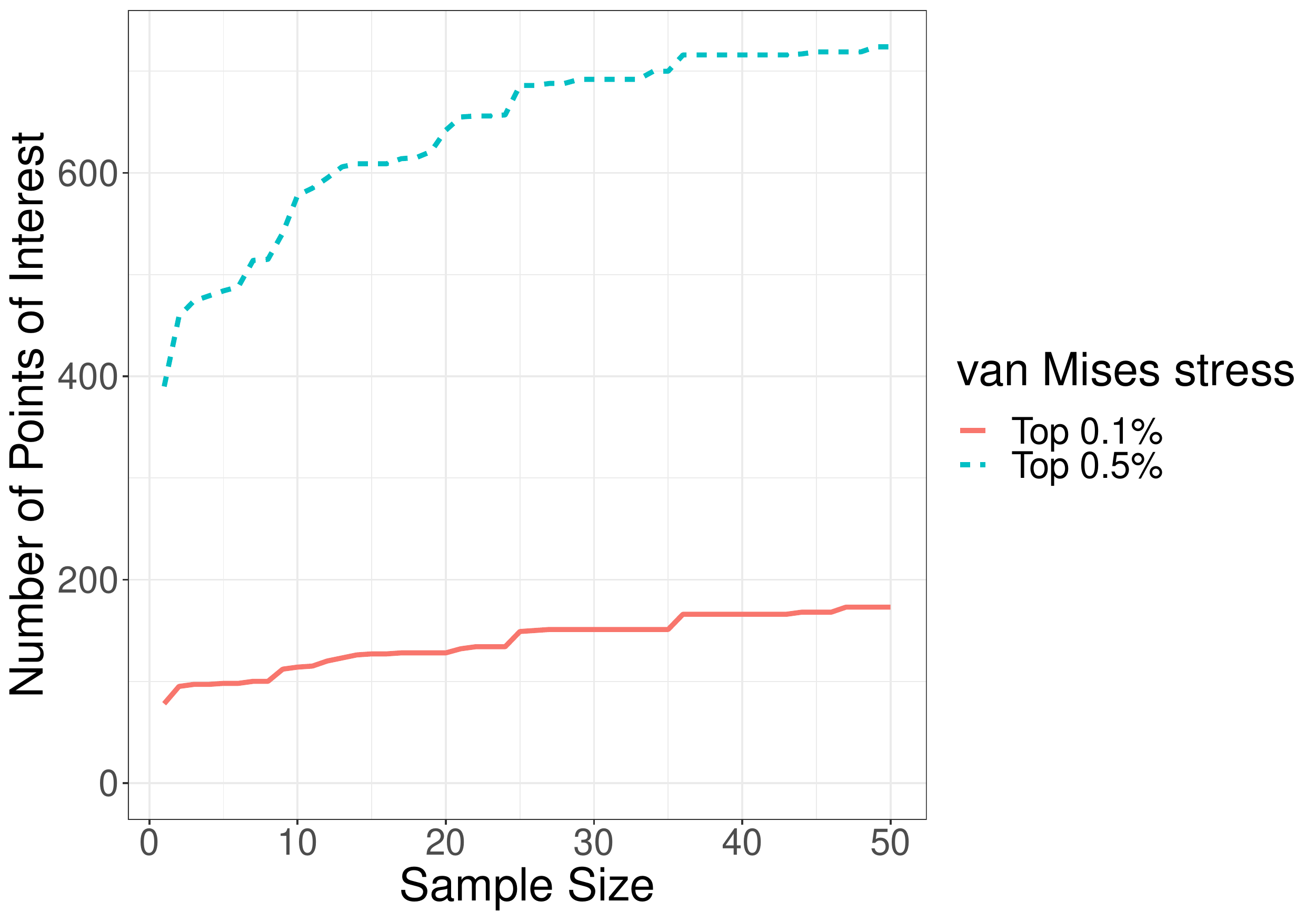}
  \end{center}
  \vspace{-10pt}
  \caption{Number of gridpoints with at least one high stress realization in the training set.}
  \label{fig: points of interest}
\end{figure}
Figure \ref{fig: points of interest}
gives the average number of points on the mesh above
the quantile when increasing the sample size.
We see that high values of van Mises stress tend to occur in similar
regions and the number of relevant points on the mesh seems to be bounded
when increasing the sample size.
This shows that the spatial variance of the high stress areas is low
and we can expect to reliably predict which parts of the grid are relevant
for high stress.

Overall, as shown in Figure \ref{fig: location high stress values},
we are dealing with locally clustered areas, i.e. some
grid points are of higher interest than others.
This is what we aim to take into account when applying dimension reduction
tools. In fact, a simple
eigenvector decomposition of the grid is less suitable to
effectively compress the relevant information.
The central contribution of this paper is therefore
to construct a suitable basis for dimension reduction
which will incorporate additional weights
to deliver local features occurring in the data at hand.

\section{Modeling of the weights}
\label{sec: Modeling of the weights}
The first step is the construction of an appropriate
weight vector. These weights should mirror the importance of
each grid point and should be chosen data driven or,
alternatively, by expert knowledge.
In the following we take a data driven approach.
The derivations are given for a general
quantile $q$.

There are multiple facts that we want to account for
when estimating the relevancy of a specific point on the grid.
Firstly, the points where high stress values occur in the training set
are most important.
Secondly, points in the local neighborhood of high stress values
are informative as well.
Finally, as we can see in the training set and as shown in Figure \ref{fig: points of interest},
high stresses do overlap between different samples and only vary slightly along
the borders. Hence, these points of the mesh are of primary interest.

Given the already high amount of information in the training data and the spatial
structure of the data we will use a smoothed version of the high stress
frequencies $\bm{f}$ as weight vector $\wbm$.
We propose standard kernel smoothing with Gaussian kernels
\begin{equation}
  K(s_{j}, s_{l})_{\sigma}
  = exp\left(-\frac{dist(s_{j},s_{i})^{2}}{2\sigma^{2}}\right).
\end{equation}
Here $dist$ defines the shortest path in edges between two points on the high fidelity grid.
The smooth weights are then defined as
\begin{equation}
  \wbm(s_{j})=\frac{\sum_{l=1}^{m_{y}}K_{\sigma}(s_{j}, s_{l})\fbm(s_{j})}{\sum_{l=1}^{m_{y}}K_{\sigma}(s_{j}, s_{l})}
\end{equation}
with $\fbm(s_{j})$ defined as in \eqref{eq: stress frequency}.
We can summarize this equation in a vector matrix form with $\wbm = \bm{S}(\sigma^{2})\fbm$.
Here $\bm{S}(\sigma^{2})$ is the smoothing matrix containing the smoothing weights.
Further we can select the bandwidth parameter $\sigma^{2}$ with standard generalized cross validation (see
\cite{Gene1979}, \cite{li1986}).
We omit the plotting of the resulting weight vectors as there is a slight
difference to Figure \ref{fig: location high stress values}.

\section{Finite Element features using additional weights}
\label{sec: Formulation of the model}
\subsection{Basis Definition}
\label{sec: Basis Definition}
Before we outline a general solution
and look at some exemplary vectors, we
define the weight based basis vectors for the high fidelity grid
as a series of optimization problems.
For this part we assume the existence of a weight vector
$\wbm \in \IR^{m_y}$ with $\sum_{i=1}^{m_{y}} \wbm(s_{j})=1$ and
$\wbm(s_{j}) \geq 0, \forall j \in \left[m_{y}\right]$. We assume further
that the high fidelity output follows a Gaussian Markov Random Field (GMRF) assumption, i.e.
$\ybm \sim N(\textbf{0}, \Sigmabm)$.
The mean vector is of secondary interest and for the covariance matrix
$\Sigmabm$ we assume
a GMRF precision matrix, $\Qbm = \Sigmabm^{-1}$, built from the finite elements structure.
$\Qbm$ has zero entries except for
\begin{align*}
  &\Qbm_{ij}= -1 \Longleftrightarrow \{s_{j}\}\ \text{in}\ ne(s_{i})\\
  &\Qbm_{ii}= |ne(s_{i})|.
\end{align*}
Here $ne(s_{i})$ defines the neighborhood of point $s_{i}$ on the grid and we refer to
\cite{Rue2005} for details.
The goal is the approximation of $\ybm$ with an orthonormal basis consisting of vectors
$\Bbm_{p} = \left(\bbm_{1},\cdots,\bbm_{p}\right)$, i.e.
$\bbm_{h} \in \IR^{m_{y}}$, with $\bbm_{h}^{t}\bbm_{h}=1$ and
$\bbm_{h}^{t}\bbm_{g}=0, \quad \text{for} \quad 1 \leq g < h \leq p$.
This basis depends on the weights and - as shown in the previous part -
the corresponding stress quantile, i.e. $\Bbm_{p}=\Bbm_{p}(\wbm(q))$.

\subsubsection{First vector}
We will compute the basis iteratively and start with a
basis of length one consisting of the single vector $\bbm$.
If we use this basis for a straight forward compression of $\ybm$ it holds
\begin{align}
  \ybmh &= \bbm\bbm^{t} \ybm \\
  \ybm - \ybmh &\sim N(0, \underbrace{(\Ibm-\bbm\bbm^{t})\Sigmabm(\Ibm-\bbm\bbm^{t})}_{=:\tilde{\Sigmabm}(\bbm)}).
\end{align}
Using this distribution we define the following squared risk function,
which summarizes the componentwise expected squared error
\begin{align}
  \label{eq: squared loss}
  \sum_{j=1}^{m_{y}} E\left(\left(\ybm(s_{j}) - \ybmh(s_{j})\right)^{2}\right) &= tr(\tilde{\Sigmabm}(\bbm)).
\end{align}
We choose the basis vector $\bbm_{1}$ to minimize
the expected L$2$ loss \eqref{eq: squared loss} and it
therefore solves the optimization problem
\begin{align}
  min_{\{\bbm\}}\left(\sum_{j=1}^{m_{y}} E\left(\left(\ybm(s_{j}) - \ybmh(s_{j})\right)^{2}\right)\right) &=
                                                                                                 min_{\{\bbm\}}\left(tr(\tilde{\Sigmabm}(\bbm))\right).
\label{eq: unweighted problem}
\end{align}
One can easily show that this problem is equivalent to an eigenvector problem, i.e. the
vector $\bbm_{1}$ that solves this optimization problem is the
eigenvector corresponding to the largest eigenvalue of $\Sigmabm$.

Instead of Formula \eqref{eq: squared loss} we now
incorporate the weight vector $\wbm$ constructed above
leading to a risk function with weighted componentwise compression errors:
\begin{align}
  \sum_{j=1}^{m_{y}} \wbm(s_{j})E\left(\left(\ybm(s_{j}) - \ybmh(s_{j})\right)^{2}\right) &= tr(\Wbm\tilde{\Sigmabm}(\bbm)),
\label{eq: weighted problem}
\end{align}
with $\Wbm = Diag(\wbm)$.
The optimal basis $\bbm_{1}$ then solves the optimization problem:
\begin{align}
  \label{eq: opt problem}
  min_{\{\bbm\}}\left(\sum_{j=1}^{m_{y}} \wbm(s_{j})E\left(\left(\ybm(s_{j}) - \ybmh(s_{j})\right)^{2}\right)\right) &=
                                                                                                          min_{\{\bbm\}}\left(tr(\Wbm\tilde{\Sigmabm}(\bbm))\right), \quad \norm{\bbm} = 1.
\end{align}
This problem is not equivalent to simple eigenvalue computation
and we have to resort to numerical optimization algorithms.
For the actual numeric computation it is preferable
to use a slightly different version of the objective risk function, which
we get from \eqref{eq: weighted problem} by plugging in $\tilde{\Sigmabm}$ and
simplifying:
\begin{align}
  \label{eq: function}
    R(\bbm):=& tr(\Wbm\Sigmabm)
    -2 \bbm^{t}\Wbm\Sigmabm\bbm
    +\bbm^{t}\Wbm\bbm\bbm^{t}\Sigmabm\bbm.
\end{align}
This leads to the gradient and the Hessian matrix
\begin{align}
  \label{eq: gradient}
  \begin{split}
    \frac{\partial R}{\partial \bbm} =& -2(\Wbm\Sigmabm + \Sigmabm\Wbm) \bbm
    + 2 (\bbm \Sigmabm \bbm) \Wbm \bbm
    + 2 (\bbm \Wbm \bbm) \Sigmabm \bbm
  \end{split}
\end{align}
\begin{align}
  \label{eq: hessian}
  \begin{split}
    \frac{\partial R}{(\partial \bbm)^{2}} =& -2(\Wbm\Sigmabm + \Sigmabm\Wbm)
    + 2 (\bbm \Sigmabm \bbm) \Wbm
    + 2 (\bbm \Wbm \bbm) \Sigmabm\\
    &+ 4 \Wbm \left(\bbm \bbm^{t}\right) \Sigmabm
    + 4 \Sigmabm \left(\bbm \bbm^{t}\right) \Wbm.
  \end{split}
\end{align}

\subsubsection{From p to (p+1)}
Given the basis vectors $\Bbm_{p} = \left(\bbm_{1} \cdots \bbm_{p}\right)$ we then derive the
basis vector $\bbm_{p+1}$ as the vector that optimally compresses the remainder and
is orthogonal to the previous basis vectors, i.e.
\begin{align}
  \begin{split}
    \ybmh_{p} &= \Bbm_{p}\Bbm_{p}^{t} \ybm \\
    \ybmh_{p+1} &= \bbm\bbm^{t} \left(\ybm - \ybmh_{p}\right).
  \end{split}
\end{align}
Then it holds
\begin{align}
\begin{split}
  \left(\ybm - \ybmh_{p}\right)- \bbm\bbm^{t}\left(\ybm - \ybmh_{p}\right)
  = \ybm - \ybmh_{p}- \bbm\bbm^{t}\ybm
  \sim N(0, \tilde{\Sigmabm}_{p}(\bbm)),\\
\end{split}
\end{align}
with $\tilde{\Sigmabm}_{p}(\bbm) =
(\Ibm-\bbm\bbm^{t})(\Ibm-\Bbm_{p}\Bbm_{p}^{t})\Sigmabm(\Ibm-\Bbm_{p}\Bbm_{p}^{t})(\Ibm-\bbm\bbm^{t})$.
Analogously to above we can define a risk function which summarizes the
weighted componentwise compression errors of the remainder $\ybm - \ybmh_{p}$
\begin{equation}
  \sum_{j=1}^{m_{y}} \wbm(s_{j})E\left(\left(\left(\ybm(s_{j})-
        \ybmh_{p}(s_{j})\right) - \ybmh_{p+1}(s_{j})\right)^{2}\right) = tr(\Wbm\tilde{\Sigmabm}_{p}(\bbm)).
\end{equation}
Given the first $p$ vectors of the basis the next vector $\bbm_{p+1}$ of the basis is then
given by an optimization problem with a different
objective function:
\begin{align}
  \label{eq: n to n+1 opt problem}
  min_{\{\bbm\}}\left(tr(\Wbm\tilde{\Sigmabm}_{p}(\bbm))\right) \quad \text{with} \quad \norm{\bbm} = 1, \quad \bbm^{t}\bbm_{h}=0 \quad
  \text{for} \quad 1\leq h \leq p.
\end{align}
If all weights are the same the optimal vector will be the $p+1$-th eigenvector of $\Sigmabm$.
A simplified version of the risk function,
the gradient and the Hessian can be computed similarly to
\eqref{eq: function}, \eqref{eq: gradient} and \eqref{eq: hessian}.
We refer to the supplied code for further details.

\subsubsection{Beyond the Gaussian}
\label{sec: beyond gaussian}
Strictly speaking, the derivations in section \ref{sec: Basis Definition}
only need the covariance matrix $\Sigmabm$. The Gaussian distribution is not
needed as all expectations in \eqref{eq: squared loss} can be expressed as linear
combinations of variances and covariances of $\ybm$.
This opens up the space for possible applications to any multivariate distribution
where the covariance matrix is known and numerically treatable.
In fact, we do not need a distribution at all and can just plug in an
empirical estimate for $\Sigmabm$.

\subsection{Heuristic Solution}
\label{sec: Heuristic Solution}
Solving the problems defined in
\eqref{eq: opt problem} and \eqref{eq: n to n+1 opt problem} is nontrivial and needs some additional attention.
First the series of optimization problems
that lead to the basis $\bbm_{1},\dots, \bbm_{p}$ are not convex.
This is also true for the
unconstrained case as well as the unweighted.
Without going into details and depending on the
nature of $\Wbm$ and $\Sigmabm$ the number of local minima is generally in the order of $2m_{y}$.
These correspond to only $m_{y}$ bases of different compression quality
due to the inherent symmetry of the problem, i.e. $R(\bbm)=R(-\bbm)$.
This means that are in general a high number of possible
local minima that render a brute force
numerical optimization approach problematic.
Standard approaches using a Lagrange multiplicator that
iterate over different points on the surface of the unit cube are
- depending on the choice of the starting point - much more likely to end up
in a local minimum of the risk function than finding the global one.
The deeper reason is the fact that
the unconstrained objective function is after all a distorted version of the
standard eigenvalue problem, i.e. the potential basis
vectors are distributed around the origin, with at most one in each direction.
The only point that gives us free access to all of them - the origin itself -
however is not part of the available set of vectors.
A possible solution is therefore to drop the normalization condition
and start numerical optimization at zero - simply following the steepest descent.
As the gradient itself vanishes at zero we will instead resort to
the smallest eigenvector of the Hessian \eqref{eq: hessian}, i.e. use
a Taylor approximation of rank two at zero.
To get a basis vector we will finally normalize the resulting vector to one.
Summing things up, this leaves us with the following heuristic for the first basis vector:

\begin{algorithmic}
  \State {1. Compute $\bm{v}_{max}$ the eigenvector with maximum eigenvalue of $\Wbm\Sigmabm + \Sigmabm\Wbm$}
  \State {2. $\bbm_{start}=\epsilon \cdot \bm{v}_{max}$}
  \State {3. Solve $min_{\{\bbm\}}\left(tr(\Wbm\tilde{\Sigmabm})\right)$ with a numerical optimizer starting at $\bbm_{start}$}
  \State {4. $\bbm_{1}=\bbm_{opt}/\norm{\bbm_{opt}}$}
\end{algorithmic}

We take a similar approach for the following optimization problems in \eqref{eq: n to n+1 opt problem}
and end up with the same algorithm but based on a different variance matrix and
including an additional step to make sure the resulting vector is orthogonal to the previously computed ones.
This algorithm is heuristic but provides - as we show in the following - good results in practice.
Note also that the broader ideas behind this approach
have been employed multiple times in the literature, see for example \cite{Cong2019}.

Finally, want to stress that analogously to the unweighted case, which is
equivalent to an eigenvector problem, our weighted basis
has a connection to tensor decomposition.
We have moved this mathematical reformulation to the Supplementary Material.

\subsection{Numerical details}
\label{sec: numerical details}
There are a number of minor numerical concerns that have to be addressed in order to make the problem solvable, which
we have moved to the Appendix.
Our heuristic itself has two steps:
In order to compute the direction of the steepest descent for the starting point,
as outlined in \ref{sec: Heuristic Solution} we use the
\texttt{RSpectra} (\cite{RSpectra}) library. This allows for the computation of
the most negative eigenvector of a matrix which is only specified by a function
defining its product with an arbitrary vector and makes it as a result possible to
compute the starting point without knowing the full Hessian.

\begin{figure}[h!]
  \begin{center}
    \includegraphics[width=0.6\textwidth]{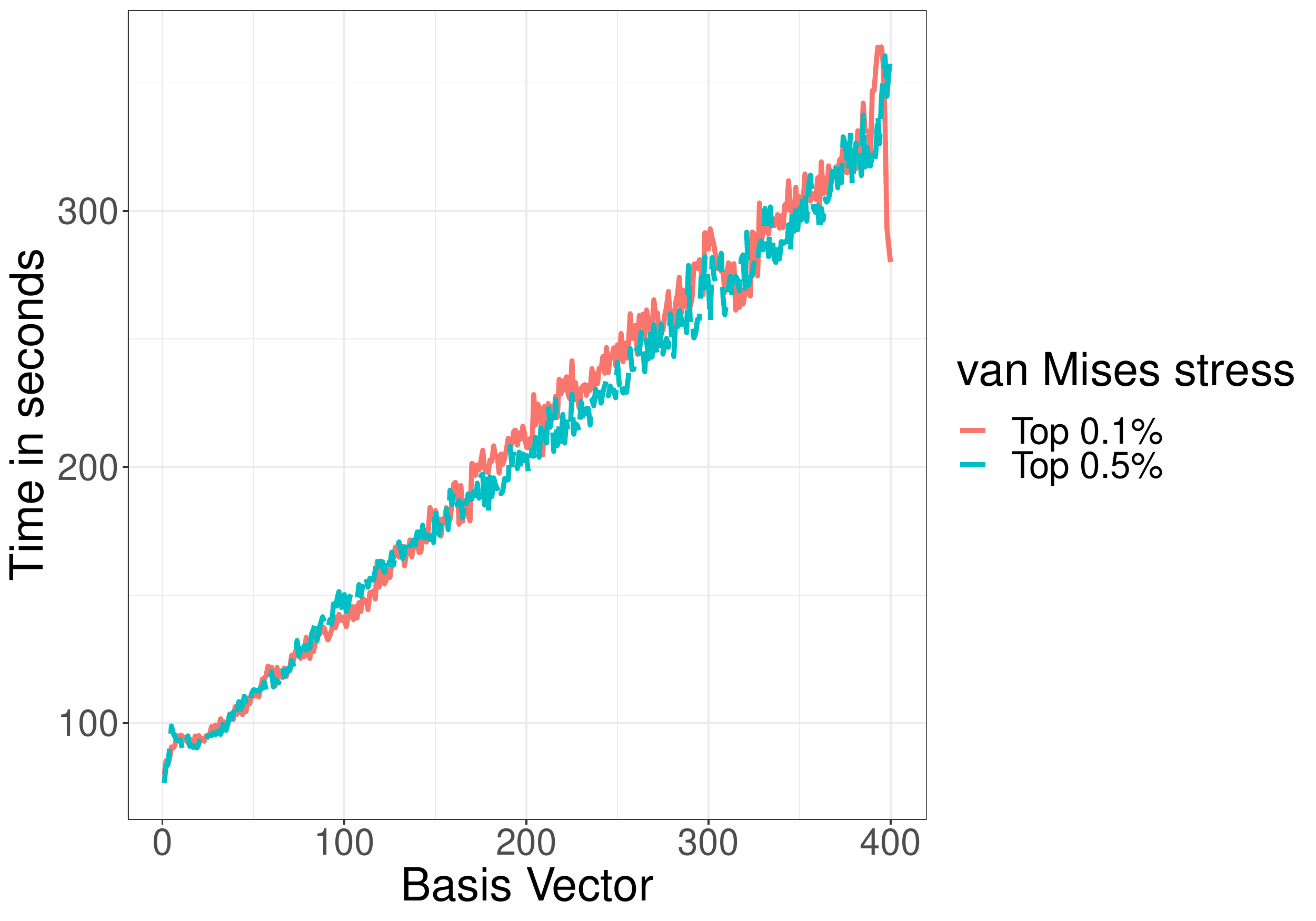}
  \end{center}
  \vspace{-10pt}
  \caption{Computation time per vector.}
  \label{fig: computation time}
\end{figure}

Given this starting point each of these problems then requires its own optimization routine.
There are a variety of solvers available for the problem at hand, which is as we have shown in
\ref{sec: Heuristic Solution} nonlinear and non-convex.
We resort to the \texttt{L-BFGS} routine which is a standard Quasi-Newton
algorithm if the gradient is available in a closed from but the Hessian isn't.
We refer to \cite{NumericalOPtimization} for a general overview of numerical optimization.
It would also be possible to use a classical Newton routine.
We plug in the weight vectors computed in Section \ref{sec: Modeling of the weights}
into our framework and
compute $400$ basis vectors for each stress quantile.
The overall computation time per vector for both weight vectors is shown in
Figure \ref{fig: computation time}.
We can see that the computation time increases linearly with the number of the computed basis vectors.
This follows directly from the problem formulation \eqref{eq: n to n+1 opt problem}: Each vector depends on all the
previously computed vectors which one by one add to the complexity of the objective function.
There is no significant difference in cost when comparing both weight vectors, even though
it might be slightly higher for the top $0.1\%$ stress values.
The overall computation time of the weighted bases is significant - on average around $16$ hours or $60000$ seconds for
a basis of length $400$.
In combination with the diminishing effectiveness of additional basis vectors this steady increase in computation time
make longer bases unattractive.

\subsection{Discussion of the Basis vectors}
\label{sec: Basis vectors}
\begin{figure}[h!]
  \begin{center}
    \includegraphics[width=0.7\linewidth]{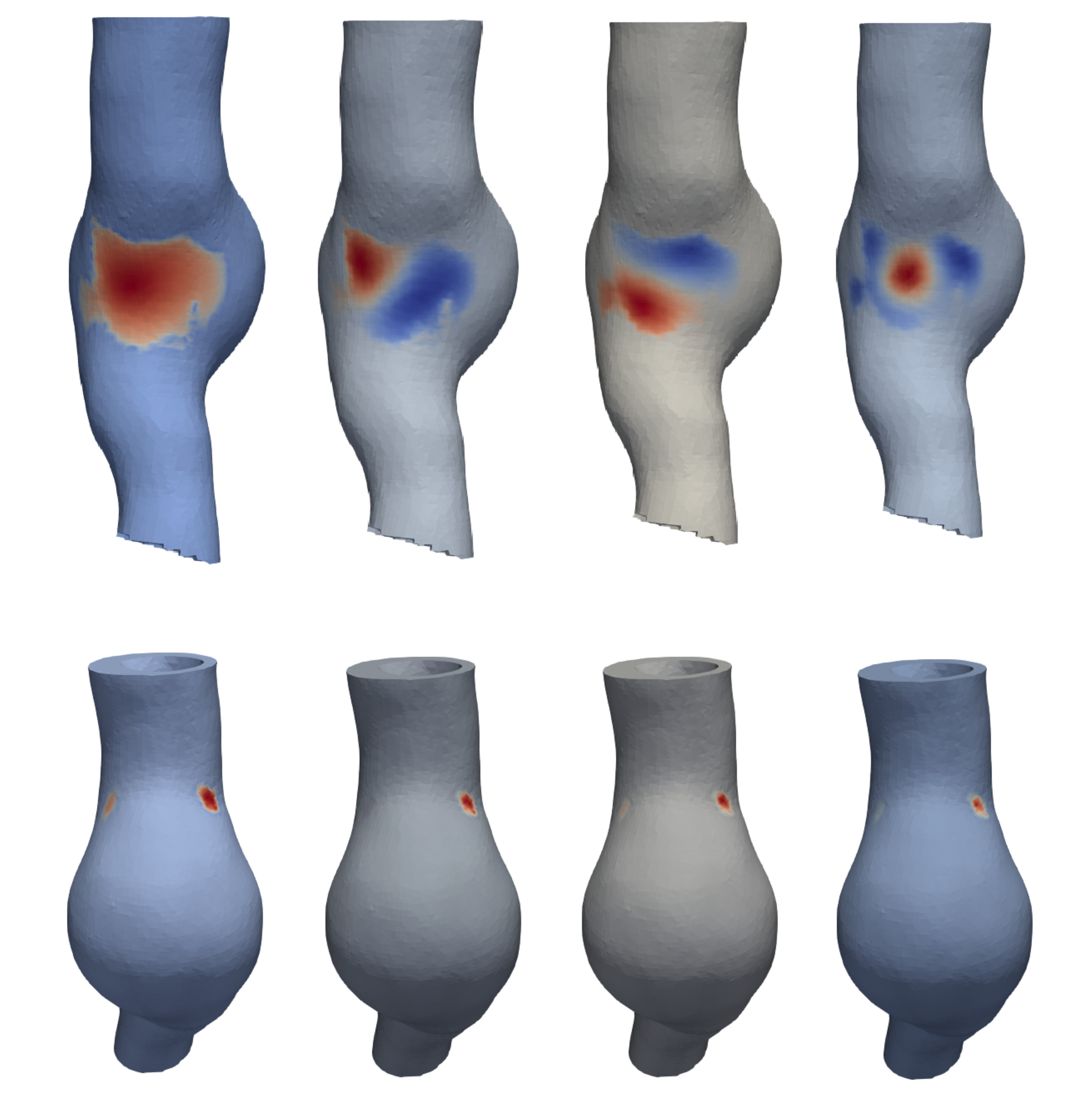}
  \end{center}
  \vspace{-10pt}
  \caption{From left to right: First four basis vectors for the weighted basis and the top $0.5\%$ stress values (frontside at the top, backside at the bottom).}
  \label{fig: basis vectors}
\end{figure}

The first four resulting vectors for the top $0.5\%$ stress values are shown in Figure \ref{fig: basis vectors}.
Basis vectors for the top $0.01\%$ values show a similar picture and are omitted.
As expected the basis vectors focus on the area in front where the
weights are bigger and then vary the pattern of the weight vector as seen as a part of
Figure \ref{fig: location high stress values}.
The remainder of the mesh grid is of little interest.

\begin{figure}[h!]
  \begin{center}
    \includegraphics[width=0.9\textwidth]{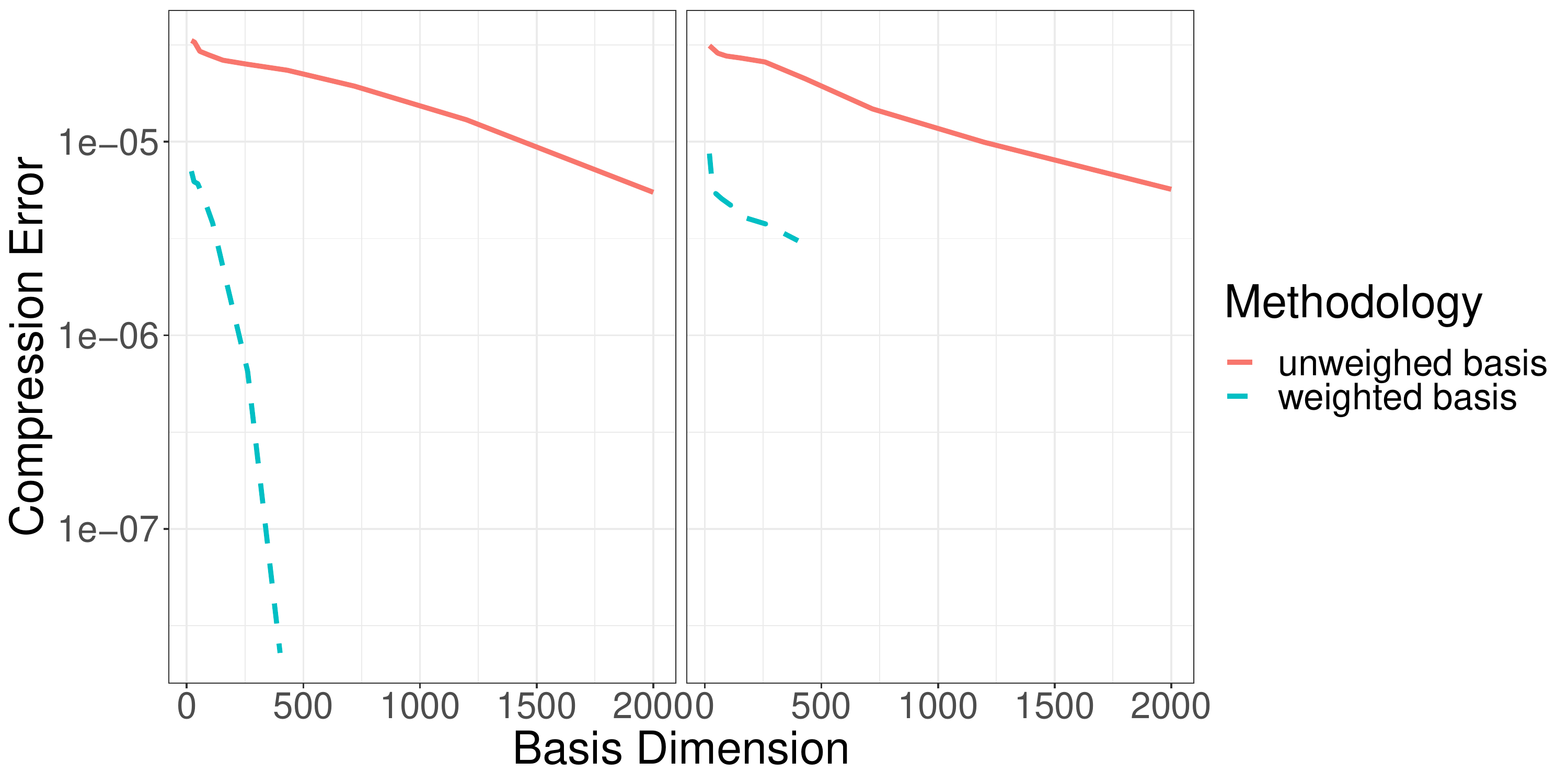}
  \end{center}
  \caption{Average squared compression error for
    $0.1\%$ (left) and $0.5\%$ (right) top stress values in the test set,
    both weighed and unweighted bases. Bases vary in length and error is shown on a
    logarithmic scale.}
  \label{fig: contained info application}
\end{figure}

We are interested in the approximation quality the bases provide, i.e. its usability for compression, especially in comparison
to ordinary unweighted GMRF bases.
These latter bases simply take the leading eigenvectors of the GMRF covariance matrix
and can be seen as a special case for a weight vector with equal weights.
We look at the empirical error
which we define as the simple mean squared compression error of the high stress values
\begin{equation}
\bm{e}_{comp}(s_{j})=\frac{\sum\limits_{k\in \mathcal{N}_{test}:
\ybm_{<k>}(s_{j})>q}(\ybm_{<k>}(s_{j})-\ybmh_{<k>}(s_{j}))^{2}}{|k\in \mathcal{N}_{test}:
\ybm_{<k>}(s_{j})>q|}.
\label{eq: compression error}
\end{equation}
Please notice that we restrict ourselves to the test set as the approximated outcomes
$\ybmh=\Bbm_{p}\left(\Bbm_{p}\right)^{t}\ybm$ depend on the weights
derived from the training set.
Further note that in practice we have to previously standardize our stress vectors by
subtracting the sample mean in correspondence with our assumptions in \ref{sec: Basis Definition}.
Results are shown in Figure \ref{fig: contained info application}.
The plot shows that the gain in compression quality decreases with each additional vector. Moreover, the
top $0.1\%$ stress values are easier to compress with a weighted basis than the top $0.5\%$.
This is to be expected:
Weighted bases get even more powerful for outcomes in a small area where the number of features
can approach the number of grid points.
For comparison we also show the compression error
for the unweighted dimension reduction, i.e. with standard eigenvectors.
Overall such comparison highlights that our
approach indeed offers a great improvement with
respect to the error. The approximation error for eigenvector bases is
an order of magnitude larger than for their weighted counterparts.
For bases of similar length we get errors which are up to three orders smaller.
The extracted eigenvectors in the unweighted basis contain
global features and the cutoff results in a
significant loss of local information.
A problem which we can solve with the introduction of weights.

\begin{figure}[h!]
  \begin{center}
    \includegraphics[width=0.9\linewidth]{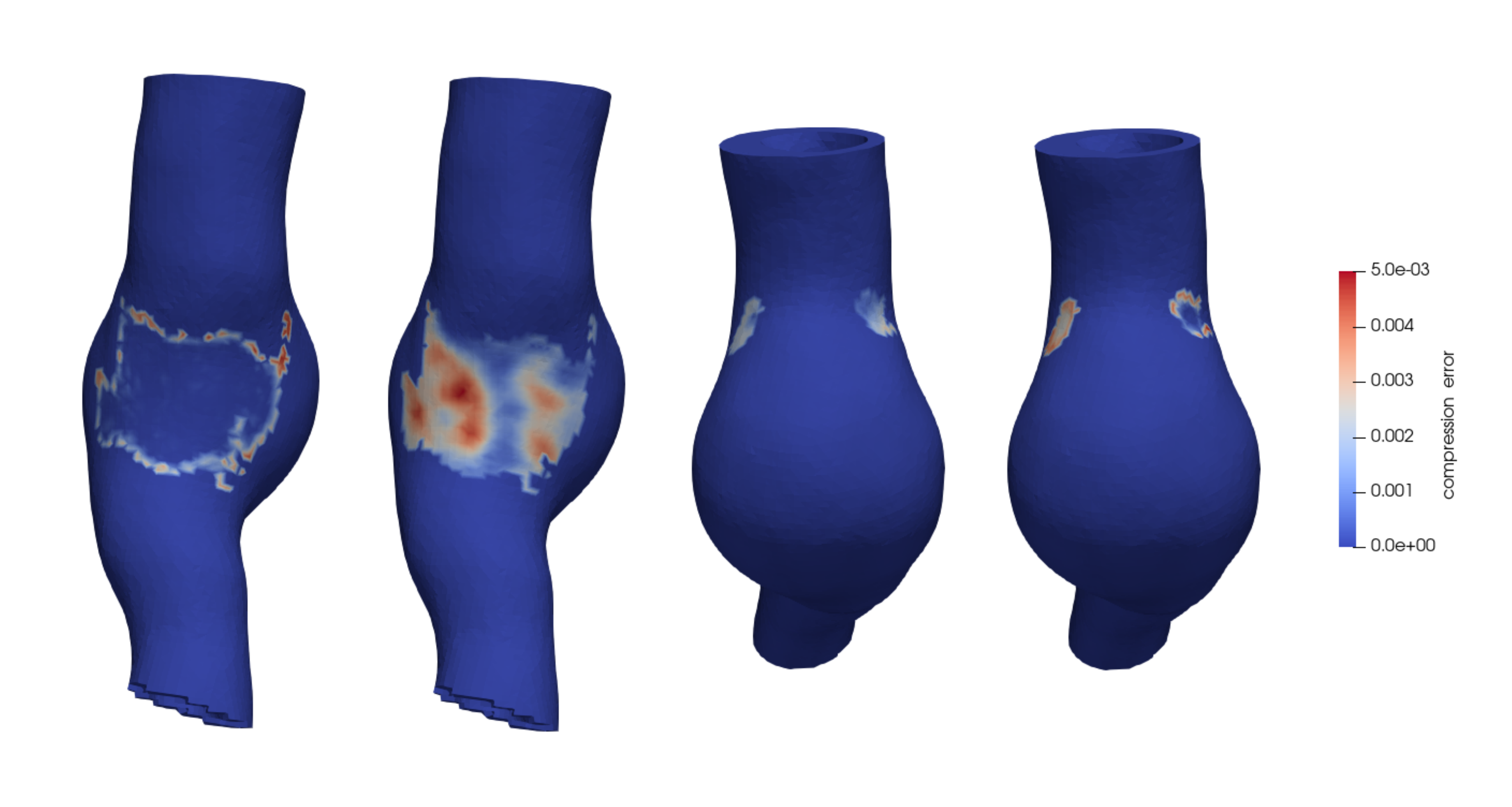}
  \end{center}
  \vspace{-10pt}
  \caption{Average squared test set compression error for weighted basis (first and third; $100$ basis vectors)
    and unweighted basis (second and fourth; $1000$ basis vectors).}
  \label{fig: 3d compression error}
\end{figure}

We also plot the empirical compression errors on top of the geometry for the top $0.5\%$ stress values
in Figure \ref{fig: 3d compression error}
to check for their spatial distribution.
In this context we remove the parts of the grid that are of no
interest for the high stress prediction, i.e. for all points
where a high stress value is not taken at least once for the $850$ test cases we set
the error to zero.
The error shows the impact of the additional weights.
While the front of the geometry has a random pattern of higher
errors for the unweighted basis,
the error decreases from the middle to the border in the case of the weighted basis.
This is explained by the weights which follow Figure \ref{fig: location high stress values} and put most emphasis in the middle.
Note that the points on the border with higher error are not necessarily problematic as the
plotted error does not adjust for the frequency of high stress values, which is low for these points.
Looking at the back of the geometry we again get a more random structure for the unweighted basis while for the weighted basis the error
is distributed in a more uniform fashion.
Overall we see that our weighted basis provides an efficient compression given enough data to narrow down
the regions on interest while not relying on enough data to estimate the actual distribution.

\section{Application: Prediction of High Stress}
\label{sec: Prediction of the maximum}
\begin{figure}[h]
  \begin{center}
    \begin{tikzpicture}[auto]
      \node (S1) {$\xbm$};
      \node (U1) [right= 2cm and 3cm of S1] {$\ubm_{x}$};
      \node (U2) [right= 2cm and 3cm of U1] {$\ubmh_{y}$};
      \node (U3) [right= 2cm and 3cm of U2] {$\tilde{\ybm}$};
      \draw[->] (S1) to node { $\left(\tilde{\Bbm}^{x}_{p_{x}}\right)^T$} (U1);
      \draw[->] (U1) to node {regression} (U2);
      \draw[->] (U2) to node {$ \tilde{\Bbm}^{y}_{p_{y}}$} (U3);
    \end{tikzpicture}
    \caption{Model Pipeline for the prediction of the complete high fidelity outcome.}
    \label{fig: tikz:pipeline}
  \end{center}
\end{figure}
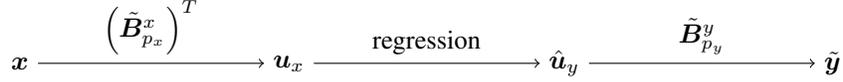

The weighted dimensionality reduction can now be utilized in practice,
where we aim to predict a high-fidelity outcome from a low-fidelity simulation.
As remarked before, when doing so we are nearly exclusively interested in large
values of the van Mises stress. We take \cite{Striegel2022} as a starting point.
That paper followed a two step approach that uses a Gaussian
Markov random field (GMRF) assumption for the construction of bases for both, low and high fidelity simulations.
These bases which we denote with $\tilde{\Bbm}^{x}_{p_{x}}$ for the $p_{x}$ dimensional low fidelity basis
and $\tilde{\Bbm}^{y}_{p_{y}}$ for its high fidelity counterpart are equivalent to the unweighted
case for the method presented in Section \ref{sec: Formulation of the model}
and consist of the eigenvectors of the GMRF covariance matrices for both grids
seen in Figure \ref{fig: low_and_high_fid_ex.pdf}.
In the second step the compressed outcomes $\ubm_{x} = \tilde{\Bbm}^{x}_{p_{x}}\left(\tilde{\Bbm}^{x}_{p_{x}}\right)^{t}\xbm$
and  $\ubm_{y} = \tilde{\Bbm}^{y}_{p_{y}}\left(\tilde{\Bbm}^{y}_{p_{y}}\right)^{t}\ybm$
were connected with a regression model as sketched in Figure \ref{fig: tikz:pipeline}.
In this paper we improve upon this approach by replacing
the high fidelity basis $\tilde{\Bbm}^{y}_{p_{y}}$ with the weighted bases $\Bbm^{y}_{p_{y}}(q)$,
which was previously denoted with $\Bbm_{p}$
in Section \ref{sec: Basis vectors} and also depends on the stress quantile $q$.

\begin{figure}[h!]
  \begin{center}
    \includegraphics[width=0.9\textwidth]{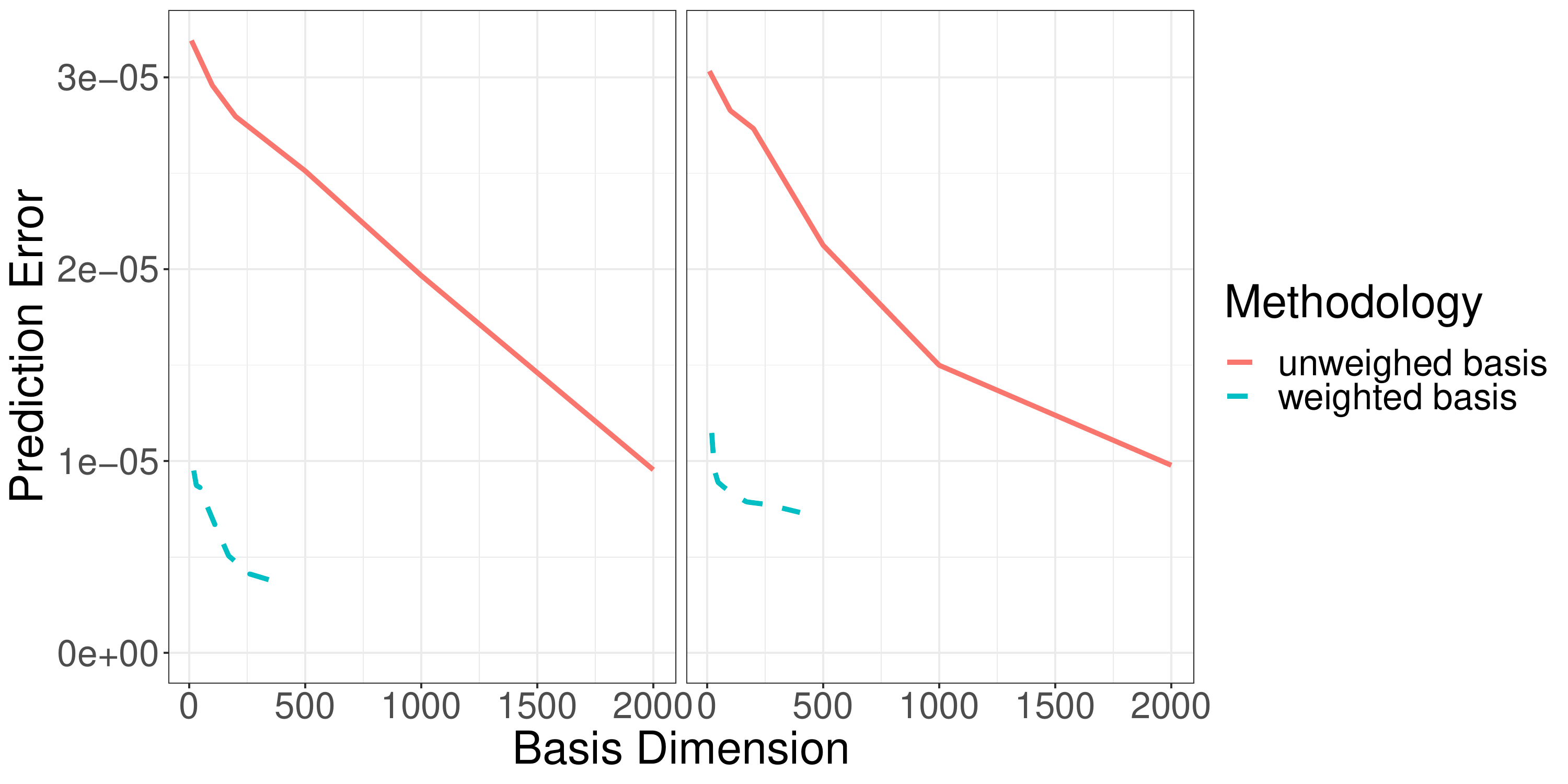}
  \end{center}
  \caption{Average squared prediction error for models fitted with different high fidelity bases lengths for
    top $0.01\%$ (left) and $0.05\%$ (right) stress values.}
  \label{fig: prediction errors}
\end{figure}

We again use a training set of $50$ cases while the remaining $850$ cases
are left for the purpose of testing. Further we use $p_x = 100$
as experiments show that these vectors are sufficient for an effective compression of
the low fidelity outcome, which is defined on an around $25$ times smaller grid.
Please see the Appendix for a plot of the low fidelity compression error versus the
(unweighted) basis length.
At the same time we vary the length of the basis used
for the high fidelity outcomes $p_{y}$.
We define the average prediction error as the mean
squared error between the top stress values in
the test set and the values predicted by our model for the relevant locations:
\begin{equation}
\bm{e}_{pred}(s_{j})=\frac{\sum\limits_{k\in \mathcal{N}_{test}:
\ybm_{<k>}(s_{j})>q}(\ybm_{<k>}(s_{j})-\tilde{\ybm}_{<k>}(s_{j}))^{2}}{|k\in \mathcal{N}_{test}:
\ybm_{<k>}(s_{j})>q|}.
\end{equation}
For comparison we compute $\bm{e}_{pred}$ for both the weighted basis $\Bbm^{y}_{p_{y}}$
and the unweighted $\tilde{\Bbm}^{y}_{p_{y}}$.
Prediction errors are shown in Figure \ref{fig: prediction errors}.
The overall shape follows the compression error shown in Figure \ref{fig: contained info application}: It is
simpler to predict the more local top $0.1\%$ values than the top $0.5\%$ given a weighted basis of fixed length.
Just as in the case of the compression the effectiveness of additional basis vectors decreases for both unweighted and
weighted bases. The absolute error however is significantly higher as it also depends on the predictive
power of the low fidelity outcomes. Finally, the advantage of the weighted basis in comparison
to its unweighted counterpart varies for different
basis lengths. As the unweighted bases approach $2000$ eigenvectors we get similar results to our weighted approach for the
top $0.05\%$ stress values. Still the weighted basis offers a slightly better performance while being four times smaller.

\begin{figure}[!]
  \centering
  \begin{subfigure}[b]{0.47\textwidth}
    \centering
    \includegraphics[width=1\textwidth]{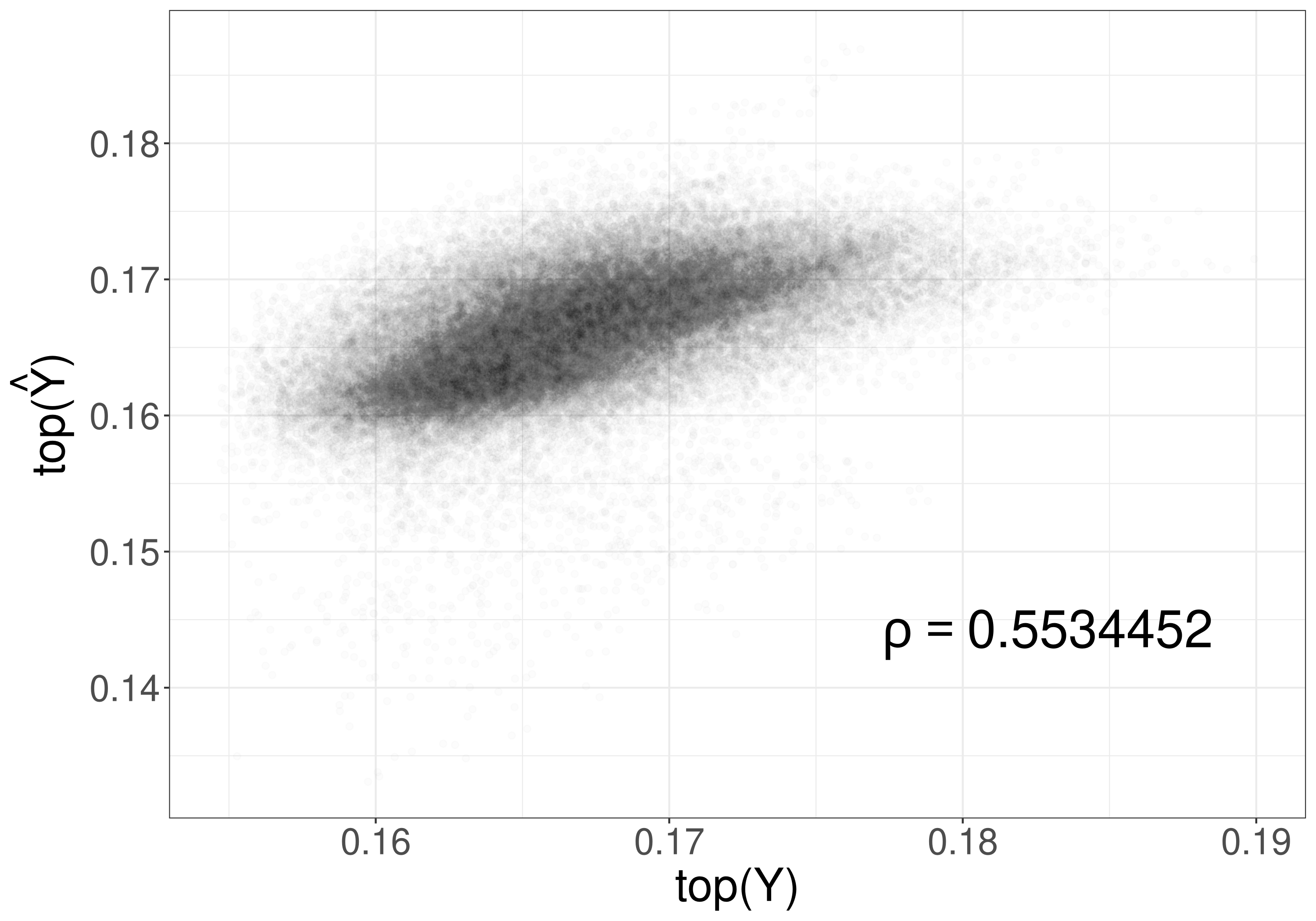}
  \end{subfigure}%
  \begin{subfigure}[b]{0.47\textwidth}
    \centering
    \includegraphics[width=1\textwidth]{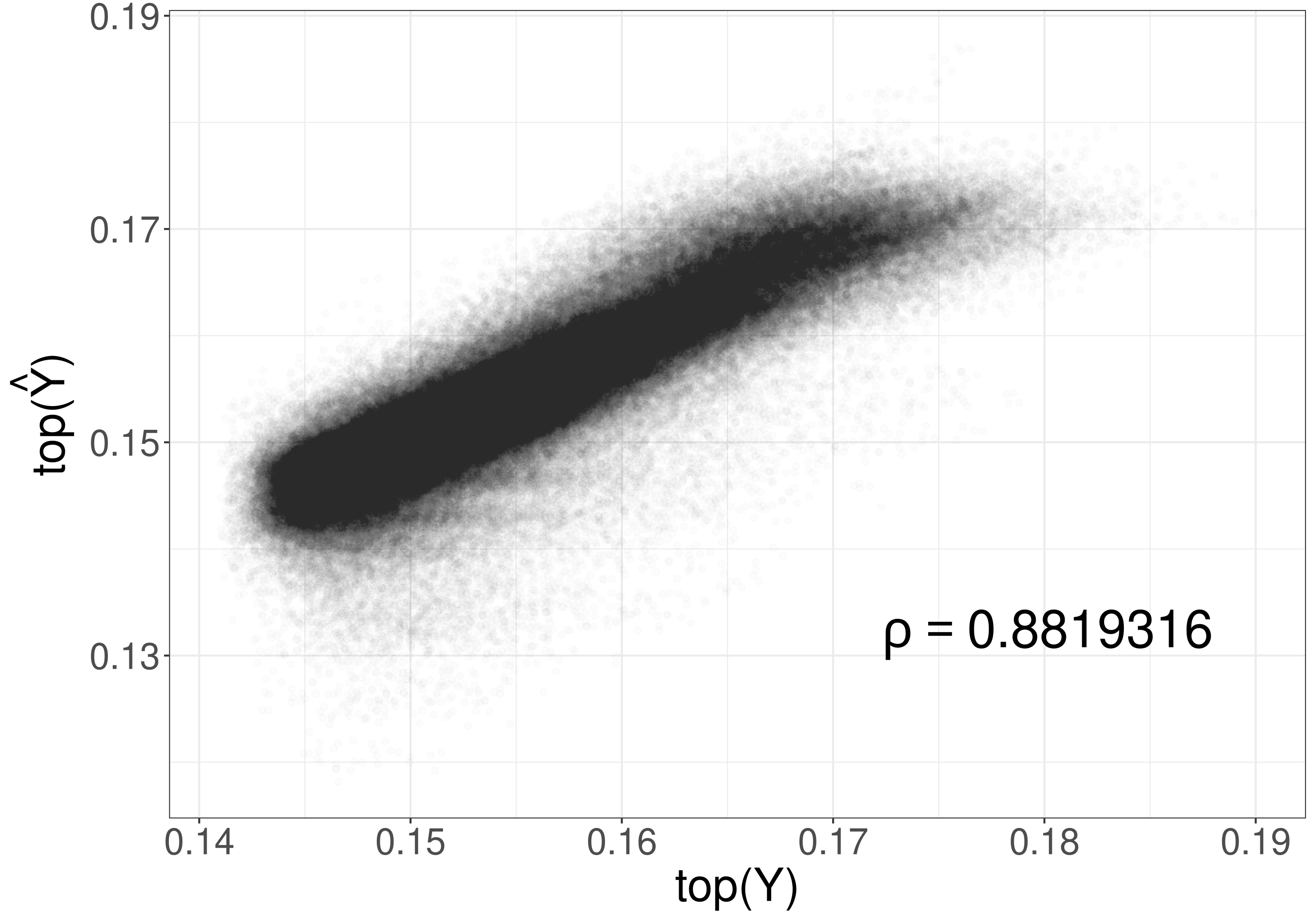}
  \end{subfigure}\\
  \begin{subfigure}[b]{0.47\textwidth}
    \centering
    \includegraphics[width=1\textwidth]{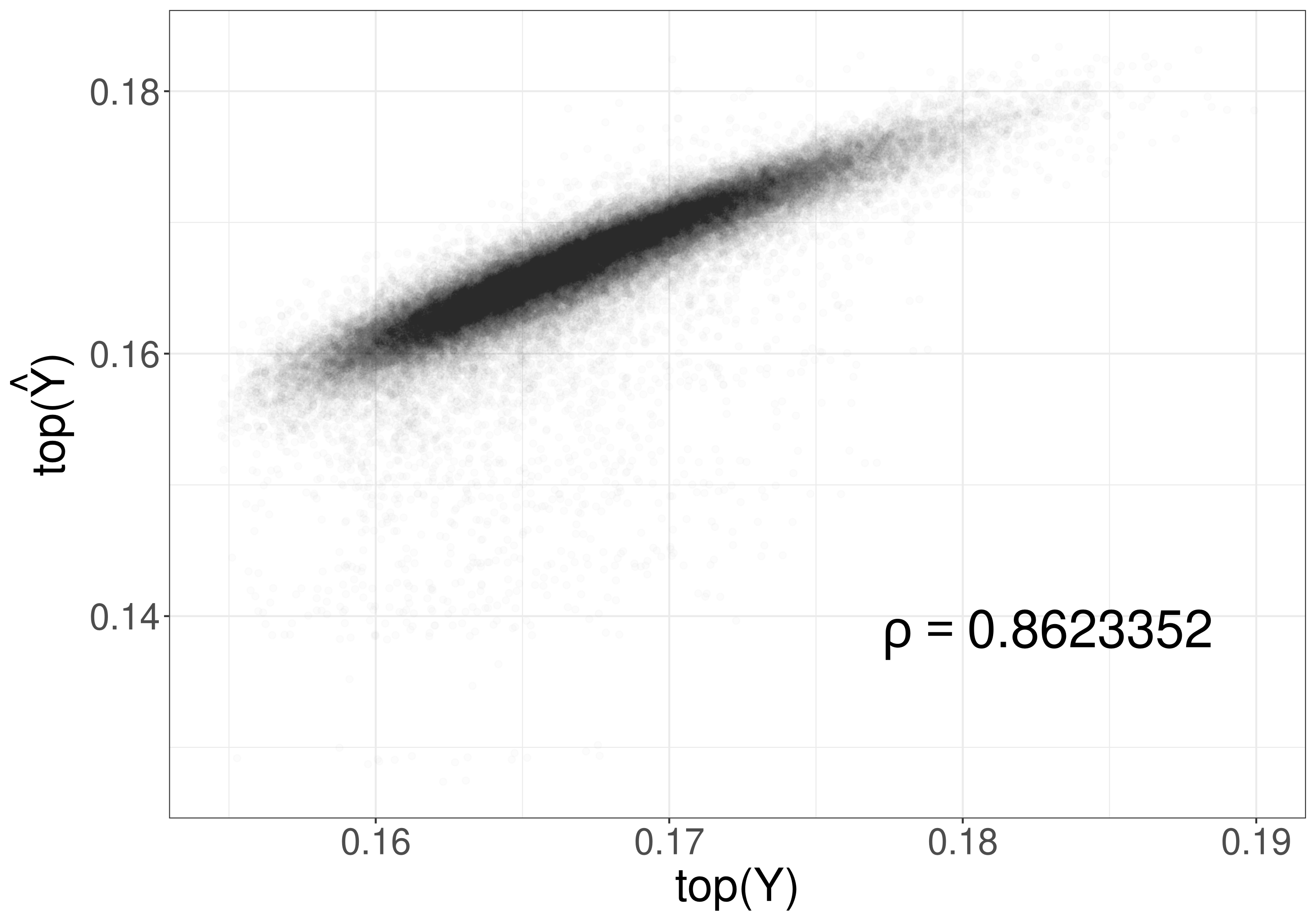}
  \end{subfigure}
  \begin{subfigure}[b]{0.47\textwidth}
    \centering
    \includegraphics[width=1\textwidth]{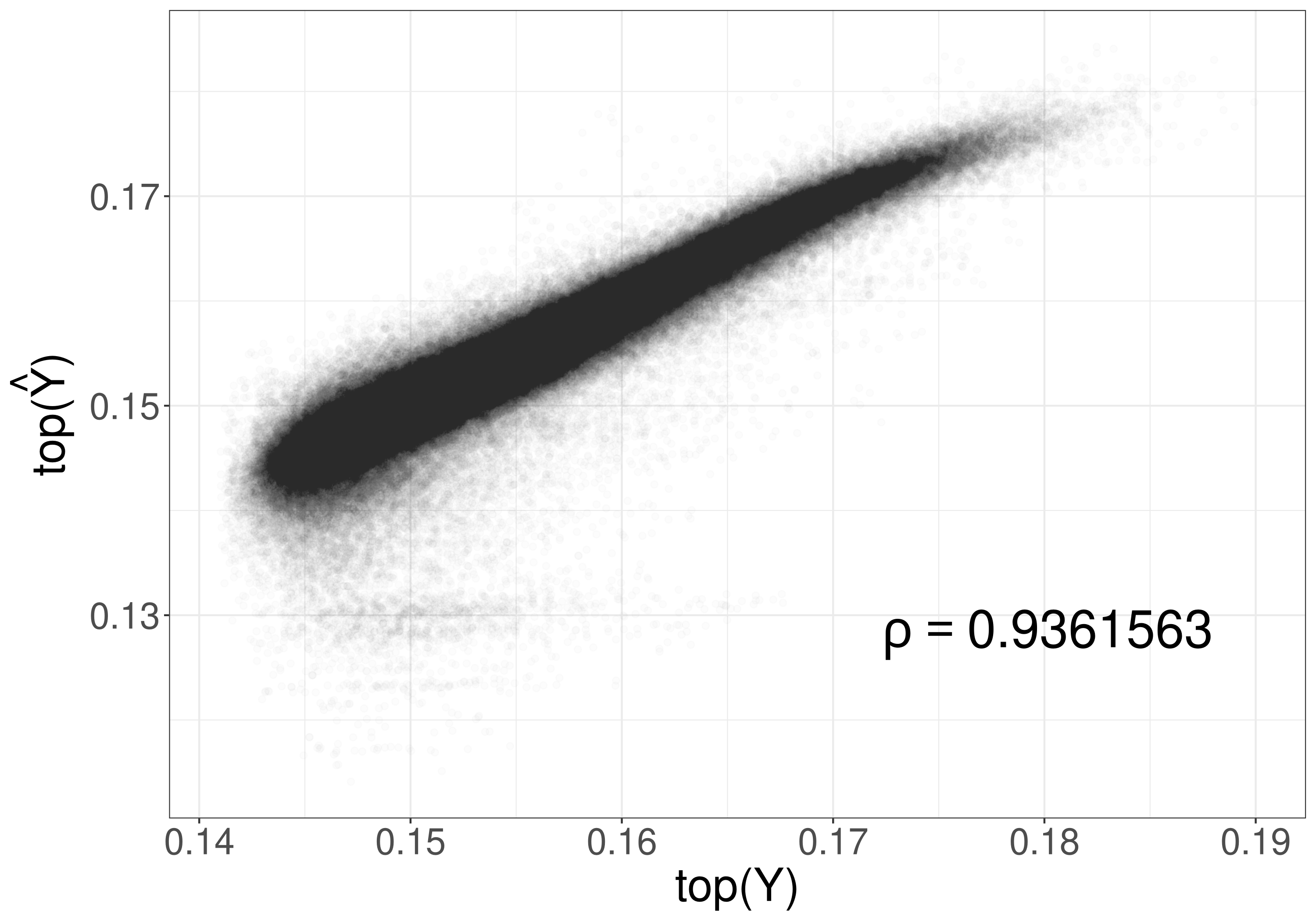}
  \end{subfigure}%
  \caption{Scatterplots of real and predicted top $0.1\%$ (leftmost) and $0.5\%$ (rightmost)
    stress values utilizing unweighted bases (top) and weighted bases (bottom) including correlation coefficients.}
  \label{fig: scatterplots}
\end{figure}

To add a further perspective on the prediction quality we also plot predicted versus real
outcomes for the weighted and unweighted bases in Figure \ref{fig: scatterplots}.
We use $p_{x}=100$ vectors for the low fidelity and
$p_{y}=1000$ for the unweighted high fidelity basis.
For comparison we add a model that only uses $p_{y}=100$ highly powerful weighted basis vectors for $\Bbm^{y}_{p_{y}}(q)$.
As we see from the two plots at the top the unweighted
approach is not well capable to predict high stress values properly - even though
there is some correlation present.
Bearing in mind that high stress levels correspond to the
risk of a highly dangerous rupture of the blood vessel
these results are not satisfying.
The main reason is the use of a basis for the high fidelity outcome
that is not suitable for the compression and prediction of the relevant local
high stress parts.
The scatterplots for the weighed bases at the bottom show that even though
there is a sizable error - most likely due to limited information
in the low fidelity outcomes - there is a significant improvement due to the use of weighed basis vectors.
This holds despite the fact that the unweighted basis is ten times longer. Points are closer to the diagonal and the
average error is lower. Nevertheless, we want to highlight that it
remains difficult to predict the highest stress levels.

\section{Discussion}
This paper proposes dimension reduction
utilizing weights corresponding to the relevancy of observations.
The starting point is a formula for the projection error using the covariance matrix,
which permits a wide range of possible applications.
We thereby take advantage of the sparse structure of its
inverse, the precision matrix - which follows from the equally sparse neighborhood structure.
This allows to numerically solve the resulting
optimization problems.
Hence, our method only requires numerically treatable
covariance structures.
Our approach shows advantages in a real world example
when predicting high fidelity outcomes from low fidelity
simulations. However, it is applicable to other settings
as well, where dimension reduction is required but data points
differ by relevance, measured by weights.

\section*{Appendix}
\subsection*{Miscellaneous Numerics}
There are a number of crucial numerical concerns that need to be
addressed in order to make the basis computation problem laid out in the previous parts
solvable in practice.
Firstly, we only know the precision matrix, $\Qbm$, but not
its inverse the GMRF covariance matrix $\Sigmabm$, as seen in \eqref{eq: function}.
$\Qbm$, the standard GMRF precision matrix (\cite{Rue2005})
is sparse.
Moreover, this matrix is not invertible in a strict sense, as the constant
vector is an eigenvector with eigenvalue $0$. We can however fix this by
adding an identity matrix, i.e. in this work we use
\begin{align}
  \tilde{\Qbm} = \Qbm + \epsilon \cdot \Ibm_{m_{y}},
\end{align}
where $\epsilon$ is a small constant, i.e. $10^{-4}$ and $\Ibm_{m_{y}}$
the $m_{y}$ dimensional identity matrix.

It is not useful nor common to actually invert $\tilde{\Qbm}$ as this is
computationally expensive and results in a dense matrix
which consumes a large amount of memory. Instead we compute
a sparse Cholesky decomposition which is enough to rapidly compute the
matrix / vector operations involving $\Sigmabm$ in the formulae for the
function value and gradient, like
\eqref{eq: function} and \eqref{eq: gradient}.
The decomposition is not sufficient to actually compute the full Hessian matrix
as seen in \eqref{eq: hessian}.
However, this is not necessary as
the ability to compute the product of the Hessian
with a given vector is sufficient in our case.

There are specialized routines for the case of GMRF type precision matrices
which do not only build on the highly sparse nature but also the (up to permutation) band
like shape of the matrix in order to achieve a computational cost of up to $O(n^{3/2})$.
This results in a large advantage with regard to computation time
in comparison to the ordinary Cholesky decomposition routines which also do not result
in sparse matrices.
We refer to \cite{DavisT} as a comprehensive reference.
For the computation we use the \texttt{R} package \texttt{Matrix}, which is based on the \texttt{CHOLMOD}
(\cite{SuiteSparse}) library.
On our machine with $64$ GB of RAM and $12$ core
AMD Ryzen $9 \quad 3900$X processor the computation of the decomposition for our $80000$ dimensional precision matrix
takes just $40$ seconds.

\subsection*{Low Fidelity Compression Error}
\begin{figure}[h!]
  \begin{center}
    \includegraphics[width=0.6\textwidth]{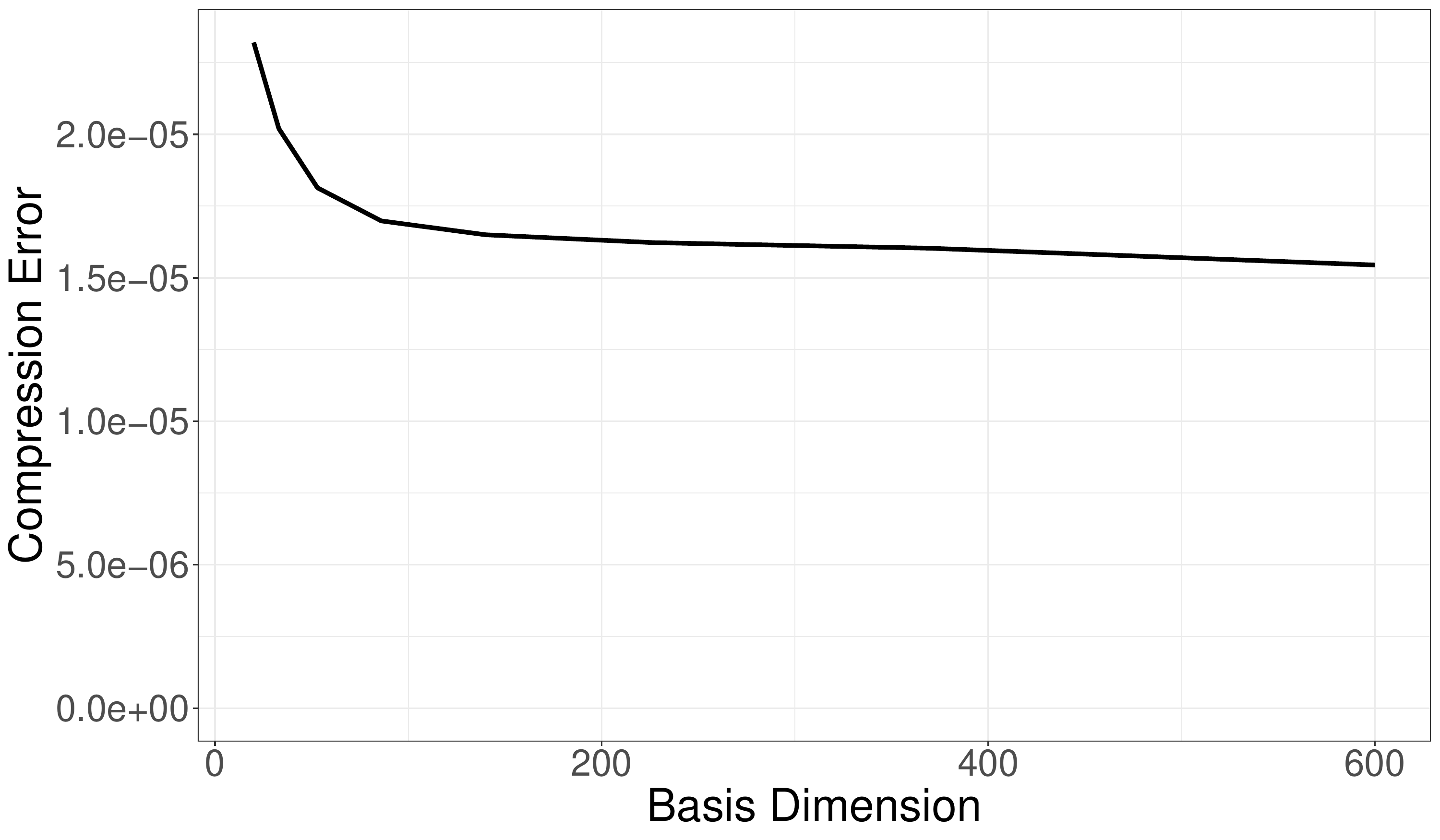}
  \end{center}
  \vspace{-10pt}
  \caption{Squared compression error for the low fidelity data using a simple (unweighted) basis of GMRF
  eigenvectors.}
  \label{fig: low fidelity compression plot}
\end{figure}

\section*{Statements and Declarations}
No funds, grants, or other support was received.
\section*{Additional Material}
\begin{description}
  \item[Code and Data:] full code and data to reproduce the results is available at:\\
        \href{https://drive.google.com/drive/folders/1c68u1bATuOZEIAXuOYzIFb3urSNQwUav?usp=sharing}{https://drive.google.com/drive/folders/1c68u1bATuOZEIAXuOYzIFb3urSNQwUav?usp=sharing}.

  \item[Reformulation as Tensor Decomposition:]  \texttt{reformulation\_tensor\_decomposition.pdf} contains
        a reformulation of the first basis vector problem as a tensor decomposition.

  \item[Simulation Study:]  \texttt{simulation\_study.pdf} contains as small simulation study with
        sampled weight vectors.
\end{description}

\bibliographystyle{spbasic}
\bibliography{./max_paper.bib}

\end{document}


\todo[inline, backgroundcolor=white]{Supplementary Material to "Weighted high dimensional data reduction of
  finite Element Features -  An Application on High Pressure of an Abdominal Aortic Aneurysm" by
Striegel et al.}
\section*{Reformulation as Tensor Decomposition}
The computation of the basis for the unweighted case
is equivalent to an eigenvector problem. Similarly the series of optimization problems discussed in
the paper can be reformulated in terms of a
tensor decomposition. We will restrict ourselves to the first basis
vector in the following as an exemplary case.
In more detail, we will show that $\bbm_{1}$ is the first vector of a
multilinear singular value decomposition for a tensor that we will derive in the following.

We define $\bbm_{1}$ as the solution of
\begin{align}
  \label{eq: opt problem}
   min_{\{\bbm\}}\left(tr(\Wbm\tilde{\Sigmabm})\right), \quad \norm{\bbm} = 1
\end{align}
with $\tilde{\Sigmabm} = (\Ibm-\bbm\bbm^{t})\Sigmabm(\Ibm-\bbm\bbm^{t}).$

Then it holds (using the simplifications from the paper):
\begin{align}
\begin{split}
  min_{\{\bbm\}}\left(tr(\Wbm\tilde{\Sigmabm})\right)
    &=min_{\{\bbm\}}\left(-2 \left(\bbm\right)^{t}\Wbm\Sigmabm\bbm
    +\left(\bbm\right)^{t}\Wbm\bbm\left(\bbm\right)^{t}\Sigmabm\bbm\right), \quad \norm{\bbm} = 1\\
  &=min_{\{\bbm\}}\left(-2 \left(\bbm\right)^{t}\bbm
    \left(\bbm\right)^{t}\Wbm\Sigmabm\bbm
    +\left(\bbm\right)^{t}\Wbm\bbm\left(\bbm\right)^{t}\Sigmabm\bbm\right), \quad \norm{\bbm} = 1.
\end{split}
\end{align}
In order to make the multilinear form more explicit one can Taylor the objective function around zero
\begin{align}
\begin{split}
  min_{\{\bbm\}}\left(\Tcal \times_{1}\bbm\times_{2}\bbm\times_{3}\bbm \times_{4}\bbm\right),
  \quad \norm{\bbm} = 1 \\
    \Tcal = \left\{ \Tcal_{ijkl}\right\} = \left\{\frac{\partial f(\bbm)}
        {\partial \bbm_{1,i} \partial \bbm_{1,j}\partial \bbm_{1,k}\partial \bbm_{1,l}}\right\}.
\label{eq: opt problem tensor}
\end{split}
\end{align}
Here $\times_{m}$ denotes the standard mode $m$ product.
The tensor $\Tcal$ is defined by
\begin{align}
\begin{split}
  \Tcal =& -4 \Ibm \times (\Wbm \Sigmabm + \Sigmabm\Wbm) -4 (\Wbm \Sigmabm + \Sigmabm\Wbm) \times \Ibm \\
         &+ 4 \Wbm \times \Sigmabm + 4 \Sigmabm \times \Wbm
           \label{eq: opt problem tensor}.
\end{split}
\end{align}

The solution of \eqref{eq: opt problem tensor}, $\bbm_{1}$, is then given by the direction $\bbm$ in which
the tensor product decreases most, i.e. for a multilinear singular value decomposition of $\Tcal$
\begin{align}
\begin{split}
  \Tcal = \Lambdabm \times_{1} \Ubm \times_{2} \Ubm \times_{3} \Ubm \times_{4} \Ubm,
\end{split}
\end{align}
where $\Lambdabm \in \IR^{m_{y}\times m_{y}\times m_{y}\times m_{y}}$ is a tensor containing
the singular values $\lambda_{1} > \lambda_{2} > ,\dots, > \lambda_{m_{y}}$
on the diagonal and $\Ubm$ is an orthonormal basis of $\IR^{m_{y}}$
we get the solution $\bbm_{1}$ as the vector in $\Ubm$ which is associated with
the most negative singular value.
We refer to \cite{Lieven_mult_sing_vale} for details on the nature of the composition.

\begin{small}
  \bibliography{max_paper.bib}
\end{small}


\todo[inline, backgroundcolor=white]{Supplementary Material to "Weighted high dimensional data reduction of
  finite Element Features -  An Application on High Pressure of an Abdominal Aortic Aneurysm" by
Striegel et al.}

\section*{Simulation Study}
In the following we present a short simulation study which explores the
utility of the presented method for basis computation in a broader setting:
We will compute bases for randomly sampled weight vectors and
evaluate the performance with the theoretical error rate defined in the
next section.

\subsection*{Theoretical Error Rate}
There is also a theoretical (compression) error ratio available for our method
which gives the theoretical remainder of the risk given a weight vector.
Such a rate can be highly useful for non-data-generated weight vectors where empirical estimates
are not available or if data for testing is scare.
If the assumptions hold this quality measure is preferable as it does not depend on the testing sample. If they
are violated we are likely to underestimate the expected error.

We derive this rate in a straight forward manner from the well-known equivalent for PCA like compression
\begin{align}
  \label{eq: error rate}
  r_{comp}=\frac{tr\left(\Wbm(\Ibm-\Bbm_{p}\left(\Bbm_{p}\right)^{t})\Sigmabm(\Ibm-\Bbm_{p}
  \left(\Bbm_{p}\right)^{t})\right)}{tr(\Wbm\Sigmabm)}.
\end{align}
The ratio is a generalization of
the PCA error rate for the case
of unequal weights, which usually consists of the quotient of
the sum of eigenvalues corresponding to the computed eigenvectors
divided by the trace of the covariance matrix. For the top $0.5\%$ values we get as an
example for $100$ vectors $\approx 0.6$, i.e. around $40$ percent of the theoretical
information is contained in this basis.

\subsection*{Error Rate Numerics}
We want to efficiently compute $tr(\Wbm\Sigmabm)$ which is the denominator
of the error rate as given in \eqref{eq: error rate} and also a constant added to the objective
function. It holds
\begin{align}
\begin{split}
  tr(\Wbm\Sigmabm) &=tr(\Wbm^{0.5}\Sigmabm\Wbm^{0.5})=tr(\Wbm^{0.5}\Qbm^{-1}\Wbm^{0.5})\\
                     &=tr(\Wbm^{0.5}\Cbm^{-1}(\Cbm)^{-T}\Wbm^{0.5})\\
                     &=tr(\Wbm^{0.5}\Cbm^{-1}(\Cbm)^{-T}\Wbm^{0.5})\\
                     &=\norm{\Wbm^{0.5}\Cbm^{-1}}_{2}^{2}.
\end{split}
\end{align}
Here $\Cbm$ denotes the lower Cholesky matrix. Unfortunately there is
no way around computing its inverse once - but this problem is solvable in a
reasonable time due to the triangular structure of the matrix.

\subsection*{Simulation Study}
Because our methodology relies on numerical heuristics
a fair assessment can only be given by actually computing bases
for different settings, i.e. computing basis vectors for
various weight vectors $\wbm$ and covariance matrices $\Sigmabm$.
In this study we will limit ourselves to the computation of bases for
a couple of weight vectors on the high fidelity grid.
In more details we will look at $8$ different weight vectors.
These weight vectors are computed in a two step approach:
First we draw a fixed number of points on the grid ($1, 3, 9, 15$).
Then we smooth the vectors that include these locations with an
eigenvector basis of $100$ or $1000$ vectors.
The goal of this approach is to explore the quality of the routine in
response to the smoothness and the distribution of the weight vector on the grid.
We compute $25$ basis vectors for each basis.

In order to give an impression of the spatial nature of the results
we plot the weight vector for $1$ center smoothed with $100$ eigenvectors and $9$ centers
smoothed with $1000$ eigenvectors along with the $3$ leading basis vectors.

\begin{figure}[h!]
  \begin{center}
    \includegraphics[width=0.8\linewidth]{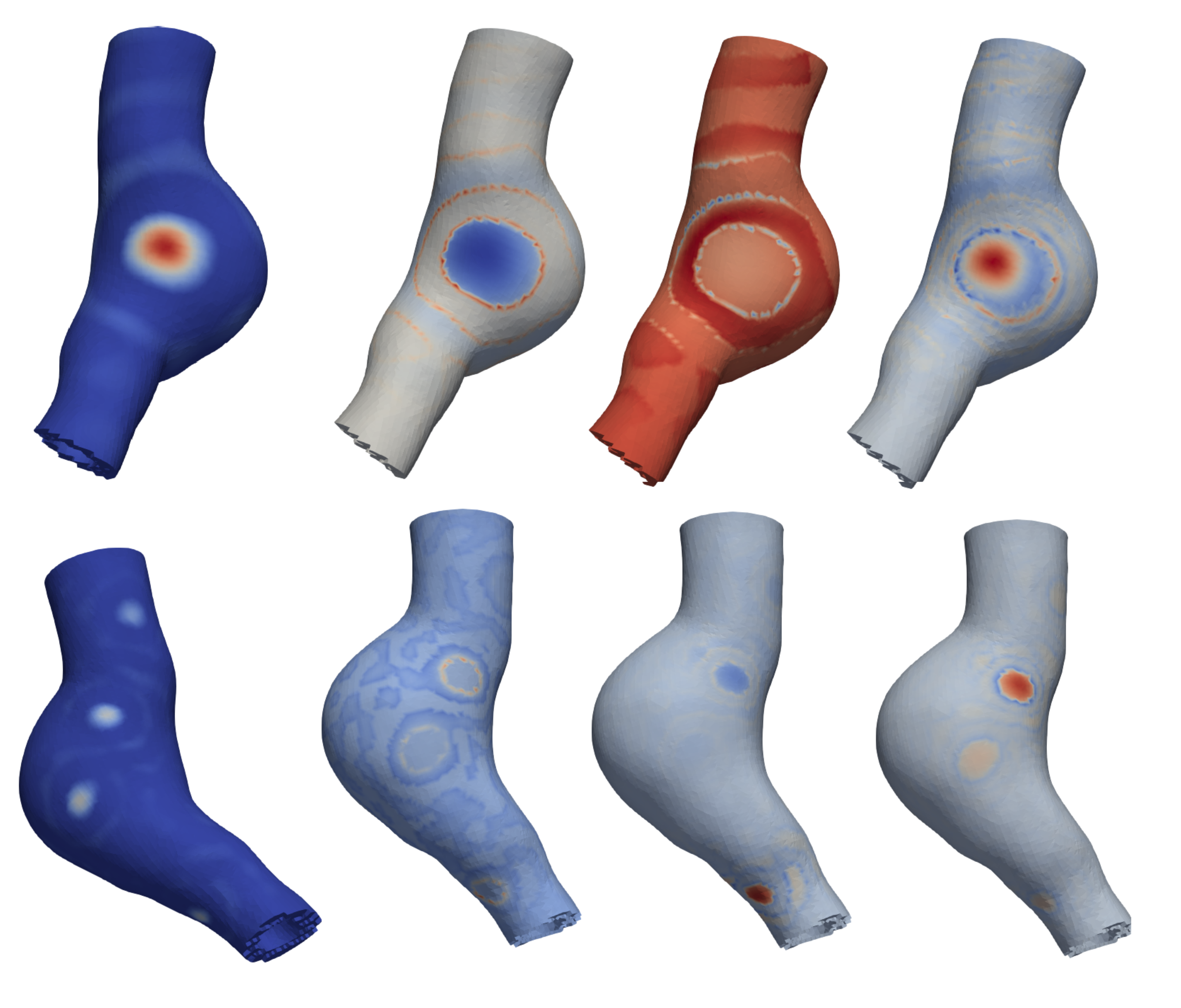}
  \end{center}
  \vspace{-10pt}
  \caption{Both top and bottom, from left to right: Weight vector
    and the corresponding first three basis vectors.
    Top we have $1$, bottom $9$ centers of weight.}
  \label{fig: observations, predicitons, errors}
\end{figure}
The basis vectors are tightly connected to the weight vectors.
They follow the pattern the weight vectors provide but alter it with each iteration.

Finally, we also take a look of the quality of the computed bases, i.e. the remainder of
\eqref{eq: error rate}, which we call the contained information $i_{c}$, versus the length of the basis.
In this context we want to point out two key facts:
First, the information rate cannot give us a proper assessment of the quality of the
numerical heuristic as we do not know the optimal basis.
Similarly, we cannot straight forward compare the information rates for different
weight vectors. The reason for this is that -
depending on the spatial distribution of the weights -
it may take more or less vectors to express the same amount of information.
In order to solve these problems one would need to adjust the information rates with
the - unknown - optimal information rate.

\begin{figure}[h!]
  \begin{center}
    \includegraphics[width=0.6\textwidth]{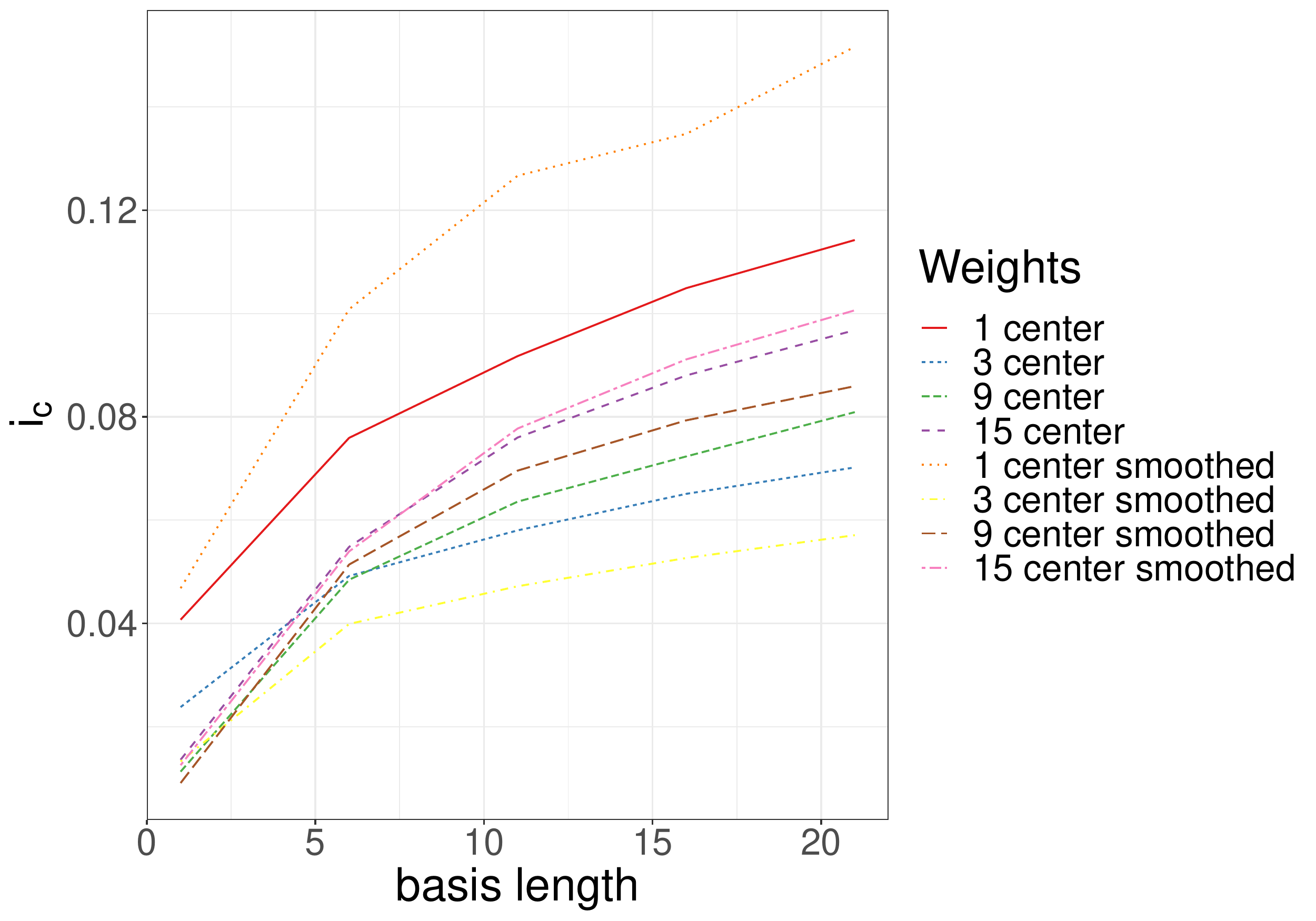}
  \end{center}
  \vspace{-10pt}
  \caption{Contained information for different bases.}
  \label{GRMF basis}
\end{figure}

Just as in the case of standard principal component analysis the improvement of the basis quality
decreases exponentially with the length of the basis.
As the first vectors contain the main features it is useful to limit oneself to
small bases in practice. It is not clear whether
more spatially distributed weight vectors make the basis computation more difficult, i.e.
our method seems to work similarly well for different weight vectors.

Even though is to be hard to interpret the absolute values the
compression rate seems satisfying given the fact that there are $80000$ grid points.
Overall the basis vectors exhibit features that are expected.

